\documentclass[a4j,11pt]{article}

\usepackage{amsmath}
\usepackage{amssymb}
\usepackage{amsthm}
\usepackage{array}
\usepackage{amscd}
\usepackage{cases}
\usepackage{bm}
\usepackage[all]{xy}
\usepackage{authblk}
\usepackage{lscape}

\theoremstyle{plain}
\newtheorem{thm}{Theorem}[section]
\newtheorem{prop}[thm]{Proposition}

\theoremstyle{definition}
\newtheorem{df}{Definition}
\newtheorem{rem}{Remark}
\newtheorem{eg}[thm]{Example}

\theoremstyle{remark}

\numberwithin{equation}{section}

\title{Four-dimensional Painlev\'e-type equations associated with ramified linear equations II: Sasano systems}
\author{Hiroshi KAWAKAMI\thanks{\texttt{kawakami@gem.aoyama.ac.jp}}}
\affil{College of Science and Engineering, Aoyama gakuin university,
5-10-1 Fuchinobe, Chuo-ku, Sagamihara-shi, Kanagawa 252-5258, Japan.}
\date{}

\begin{document}
\maketitle

\begin{abstract}
This is a continuation of the paper
``Four-dimensional Painlev\'e-type equations associated with ramified linear equations I: Matrix Painlev\'e systems'' \cite{K1}.
In this series of three papers we aim to construct the complete degeneration scheme of
four-dimensional Painlev\'e-type equations.
In the present paper, we construct the degeneration scheme of what we call the Sasano system.
\end{abstract}


\tableofcontents

\section{Introduction}\label{sec:intro}
This is the second part of a series of three papers on four-dimensional Painlev\'e-type equations
associated with ramified linear equations.
Let us outline the background of the study.


The Painlev\'e equations were found by Painlev\'e~\cite{P} and Gambier~\cite{Gm} at the beginning of the twentieth century,
as a means to obtain new special functions.
The Painlev\'e equations were originally classified into six equations.
However, the geometry of the space of initial values~\cite{O1} shows that it is natural to classify them into eight equations \cite{Sak1}.

The Painlev\'e equations can be written in Hamiltonian form (\cite{Ok, OKSO}).
Here we list the eight Hamiltonians instead of equations themselves:
{\allowdisplaybreaks
\begin{align*} 
& t(t-1)H_\mathrm{VI}\left({\alpha , \beta\atop\gamma, \delta};t;q,p\right)\\
&\quad=q(q-1)(q-t)p^2+\{ \delta q(q-1)-(2\alpha +\beta +\gamma +\delta )q(q-t)+\gamma
 (q-1)(q-t)\} p\\
&\hspace{12cm}+\alpha (\alpha +\beta )(q-t),\\ 
& tH_\mathrm{V}\left({\alpha , \beta \atop \gamma };t;q,p\right)=p(p+t)q(q-1) 
+\beta pq+\gamma p-(\alpha +\gamma )tq,\\ 
& H_\mathrm{IV}\left(\alpha , \beta;t;q,p\right)=pq(p-q-t)+\beta p+\alpha q,\quad
tH_{\mathrm{III}(D_6)}\left(\alpha , \beta ;t;q,p\right)=p^2q^2-(q^2-\beta q-t)p-\alpha q, \\
&tH_{\mathrm{III}(D_7)}\left(\alpha;t;q,p\right)=p^2q^2+\alpha qp+tp+q,\quad
tH_{\mathrm{III}(D_8)}\left(t;q,p\right)=p^2q^2+qp-q-tq^{-1}, \\ 
& H_\mathrm{II}\left(\alpha;t;q,p\right)=p^2-(q^2+t)p-\alpha q, \quad \quad 
H_\mathrm{I}\left(t;q,p\right)=p^2-q^3-tq,
\end{align*}
}
where $H_{\mathrm{J}}$ is the Hamiltonian corresponding to the Jth Painlev\'e equation.

Around the same time as the discovery of the Painlev\'e equations,
R. Fuchs~\cite{F} and Garnier~\cite{G} showed that the Painlev\'e equations can be derived from \textit{isomondromic deformations} of certain
second order linear differential equations.
This connection of the Painlev\'e equations with linear differential equations are crucially important.

Among the Painlev\'e equations, the sixth Painlev\'e equation is the master equation in the sense that
one can obtain all other Painlev\'e equations from the sixth Painlev\'e equation through degeneration.
The relationships among the Painlev\'e equations through degeneration are shown in the scheme below.
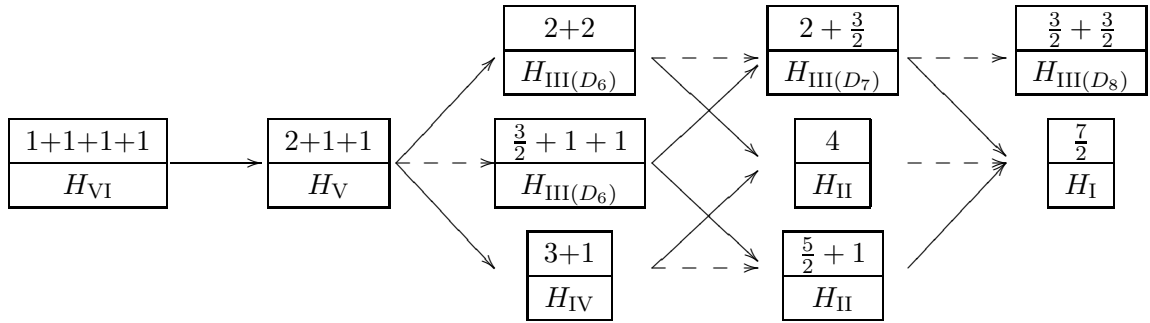
\vspace{3mm}


%

\begin{xy}
{(3,0) *{\begin{tabular}{|c|}
\hline
1+1+1+1\\
\hline
$H_{\mathrm{VI}}$\\
\hline
\end{tabular}
}},
{(35,0) *{\begin{tabular}{|c|}
\hline
2+1+1\\
\hline
$H_{\mathrm{V}}$\\
\hline
\end{tabular}
}},
{\ar (14,0);(26,0)},
{\ar (44,0);(57,14)},
{\ar @{-->} (44,0);(57,0)},
{\ar (44,0);(57,-14)},
{(67,15) *{\begin{tabular}{|c|}
\hline
2+2\\
\hline
$H_{\mathrm{III}(D_6)}$\\
\hline
\end{tabular}}},
{(67,0) *{\begin{tabular}{|c|}
\hline
$\frac{3}{2}+1+1$\\
\hline
$H_{\mathrm{III}(D_6)}$\\
\hline
\end{tabular}}},
{\ar @{-->} (78,14);(92,14)},
{\ar (78,14);(92,1)},
{\ar (78,0);(92,13)},
{\ar (78,0);(92,-13)},
{\ar (78,-14);(92,-1)},
{\ar @{-->} (78,-14);(92,-14)},
{(67,-15) *{\begin{tabular}{|c|}
\hline
3+1\\
\hline
$H_{\mathrm{IV}}$\\
\hline
\end{tabular}}},
{(102,15) *{\begin{tabular}{|c|}
\hline
$2+\frac{3}{2}$\\
\hline
$H_{\mathrm{III}(D_7)}$\\
\hline
\end{tabular}}},
{(102,0) *{\begin{tabular}{|c|}
\hline
4\\
\hline
$H_{\mathrm{II}}$\\
\hline
\end{tabular}}},
{\ar (112,14);(125,1)},
{\ar @{-->} (112,14);(125,14)},
{\ar @{-->} (112,0);(125,0)},
{\ar (112,-14);(125,0)},
{(102,-15) *{\begin{tabular}{|c|}
\hline
$\frac{5}{2}+1$\\
\hline
$H_{\mathrm{II}}$\\
\hline
\end{tabular}}},
{(135,15) *{\begin{tabular}{|c|}
\hline
$\frac{3}{2}+\frac{3}{2}$\\
\hline
$H_{\mathrm{III}(D_8)}$\\
\hline
\end{tabular}}},
{(135,0) *{\begin{tabular}{|c|}
\hline
$\frac{7}{2}$\\
\hline
$H_{\mathrm{I}}$\\
\hline
\end{tabular}}},
\end{xy}

\vspace{3mm}

The number placed in the upper half of each box is the \textit{singularity pattern}
of the associated linear equation: 
each number partitioned by + expresses the number ``Poincar\'e rank +1'' of each singular point of a linear differential equation.
The definition of the Poincar\'e rank is given in Section~\ref{sec:HTLform}.

The two kinds of arrows express the two kinds of degenerations for linear differential equations.
Namely, the dashed arrows in the above scheme correspond to degeneration of HTL canonical forms
(where HTL is an abbreviation for Hukuhara-Turrittin-Levelt, see Definition~\ref{df:HTL}) at irregular singular points,
while the others correspond to confluence of singular points of associated linear equations.


What is important is that one can obtain all the Painlev\'e equations 
by considering confluence of singular points and degeneration of HTL canonical forms of associated linear equations~\cite{OO}.



Recently, many generalizations of Painlev\'e equations based on various aspects
have been proposed and studied.
The point here is that they can be written as compatibility conditions of
some linear differential equations.
Namely, they can be regarded as isomonodromic deformation equations
(or we use the term \textit{Painlev\'e-type equations} synonymously).
Our purpose is to understand them using associated linear equations.

The scheme that we have seen in the above is the degeneration scheme of the isomonodromic deformation equations
whose phase spaces are two-dimensional.
Now we are interested in the degeneration scheme of
the Painlev\'e-type equations of four-dimensional phase space,
which comes next to the two-dimensional case.

In \cite{KNS}, the authors construct the degeneration scheme of four-dimensional Painlev\'e-type equations
by means of confluence of singular points of associated linear equations.
As the result, they obtained the degeneration scheme of four-dimensional Painlev\'e-type equations
associated with \textit{unramified} linear differential equations.

The aim of this series of papers is, by considering degeneration of HTL forms,
to obtain the degeneration scheme of four-dimensional Painlev\'e-type equations that includes ones
associated with ramified linear equations.
This gives a unified perspective of four-dimensional analogues of the Painlev\'e equations.

There are four Painlev\'e-type equations which should be placed at the starting point of
the degeneration scheme of four-dimensional case;
that is, the Garnier system in two variables~\cite{G}, the Fuji-Suzuki system (FS system for short)~\cite{FS1,Ts}, the Sasano system~\cite{Ss},
and the sixth matrix Painlev\'e system~\cite{B2,K}.
The degeneration scheme of the sixth matrix Painlev\'e system is given in our previous paper~\cite{K1}.
In the present paper, we construct the degeneration scheme of the Sasano system.

This paper is organized as follows.
In Section~2, we briefly recall HTL canonical forms of linear differential equations.
We also see the examples of a shearing transformation and make some remarks on degeneration of HTL forms.
In Section~3 we give the ramified Lax pairs and associated Hamiltonians which can be derived from
the Sasano system.
Section~4 is devoted to the description of correspondences of linear equations through the Laplace transform.

We give here the degeneration scheme of the Sasano system,
which is the main result of the present paper.

\vspace{2mm}
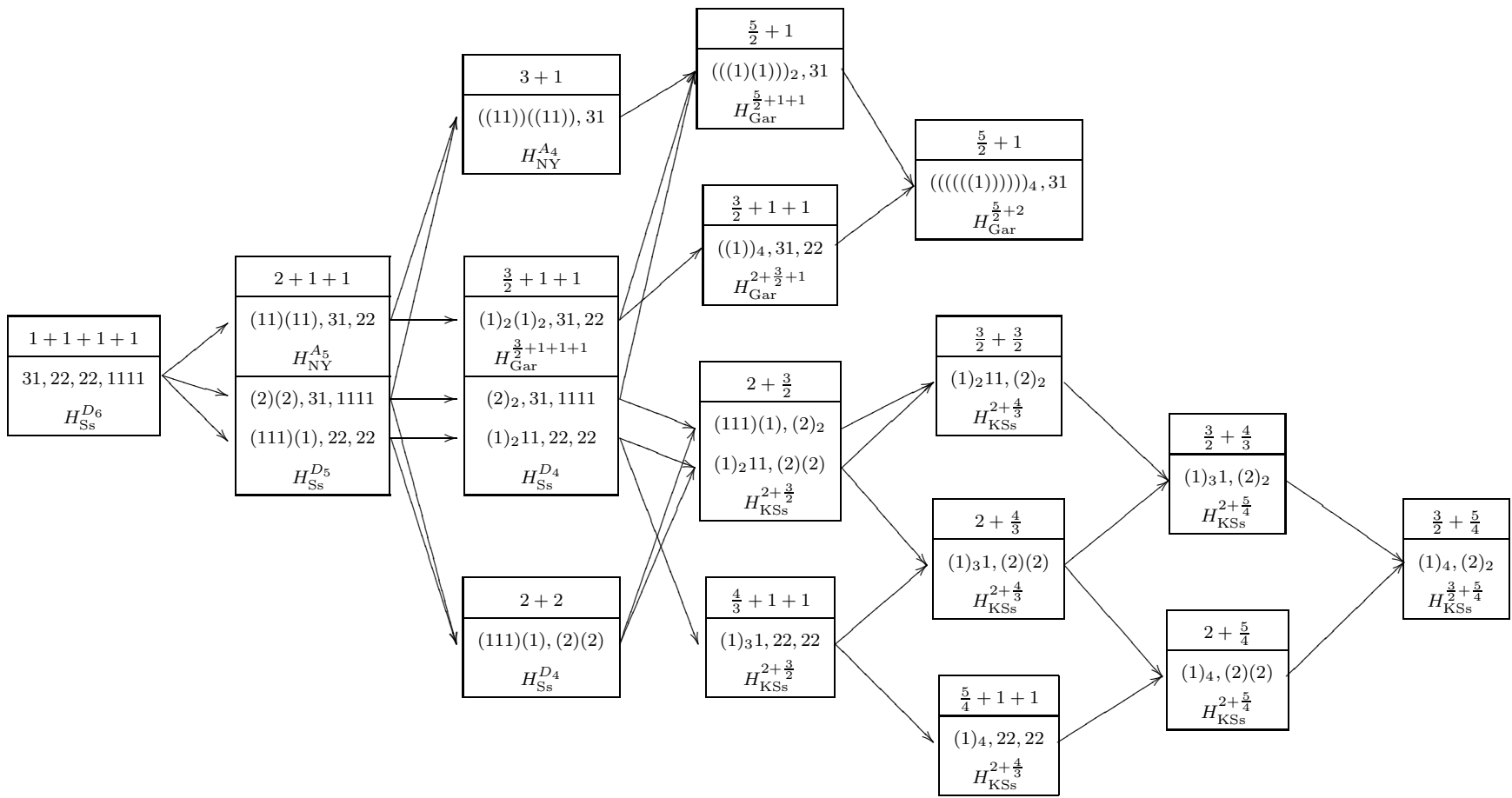
\begin{landscape}
{\scriptsize
\begingroup
\renewcommand{\arraystretch}{1.5}
\begin{xy}
{(0,-10) *{\begin{tabular}{|c|}
\hline
$1+1+1+1$\\
\hline
$31,22,22,1111$\\
$H_{\mathrm{Ss}}^{D_6}$\\
\hline
\end{tabular}}
},
{\ar (12,-10);(22,-2)},
{\ar (12,-10);(22,-13)},
{\ar (12,-10);(22,-20)},
{(35,-10) *{\begin{tabular}{|c|}
\hline
$2+1+1$\\
\hline
$(11)(11),31,22$\\
$H_{\mathrm{NY}}^{A_5}$\\ \hline
$(2)(2),31,1111$\\
$(111)(1),22,22$\\
$H_{\mathrm{Ss}}^{D_5}$\\
\hline
\end{tabular}}
},
{\ar (47,-1.5);(57,29.5)},
{\ar (47,-1.5);(57,-1.5)},
{\ar (47,-13.5);(57,29.5)},
{\ar (47,-13.5);(57,-13.5)},
{\ar (47,-13.5);(57,-51)},
{\ar (47,-19.5);(57,-19.5)},
{\ar (47,-19.5);(57,-51)},
{(70,30) *{\begin{tabular}{|c|}
\hline
$3+1$\\
\hline
$((11))((11)),31$\\
$H_{\mathrm{NY}}^{A_4}$\\
\hline
\end{tabular}}
},
{\ar (82,29.5);(93.5,36.5)},
{(70,-50) *{\begin{tabular}{|c|}
\hline
$2+2$\\
\hline
$(111)(1), (2)(2)$\\
$H_{\mathrm{Ss}}^{D_4}$\\
\hline
\end{tabular}}
},
{\ar (82,-51);(93.5,-24)},
{\ar (82,-51);(93.5,-18)},
{(70,-10) *{\begin{tabular}{|c|}
\hline
$\frac32+1+1$\\
\hline
$(1)_2(1)_2, 31, 22$\\
$H_{\mathrm{Gar}}^{\frac32+1+1+1}$\\
\hline
$(2)_2, 31, 1111$\\
$(1)_2 11, 22, 22$\\
$H_{\mathrm{Ss}}^{D_4}$\\
\hline
\end{tabular}}
},
{\ar (82,-1.5);(94.5,9.5)},
{\ar (82,-1.5);(93.5,36.5)},
{\ar (82,-13.5);(93.5,-18)},
{\ar (82,-13.5);(93.5,36.5)},
{\ar (82,-19.5);(93.5,-24)},
{\ar (82,-19.5);(94,-51)},
{(105,37) *{\begin{tabular}{|c|}
\hline
$\frac52+1$\\
\hline
$(((1)(1)))_2, 31$\\
$H_{\mathrm{Gar}}^{\frac52+1+1}$\\
\hline
\end{tabular}
}},
{\ar (116,37);(127,19)},
{(105,10) *{\begin{tabular}{|c|}
\hline
$\frac32+1+1$\\
\hline
$((1))_4, 31, 22$\\
$H_{\mathrm{Gar}}^{2+\frac32+1}$\\
\hline
\end{tabular}
}},
{\ar (115,10);(127,19)},
{(105,-20) *{\begin{tabular}{|c|}
\hline
$2+\frac32$\\
\hline
$(111)(1), (2)_2$\\
$(1)_2 11, (2)(2)$\\
$H_\mathrm{KSs}^{2+\frac32}$\\
\hline
\end{tabular}}
},
{\ar (116,-18);(130,-11)},
{\ar (116,-24);(130,-11)},
{\ar (116,-24);(129,-39)},
{(105,-50) *{\begin{tabular}{|c|}
\hline
$\frac43+1+1$\\
\hline
$(1)_3 1, 22, 22$\\
$H_\mathrm{KSs}^{2+\frac32}$\\
\hline
\end{tabular}}
},
{\ar (115,-51);(129,-39)},
{\ar (115,-51);(130,-66)},
{(140,20) *{\begin{tabular}{|c|}
\hline
$\frac52+1$\\
\hline
$((((((1))))))_4, 31$\\
$H_\mathrm{Gar}^{\frac52+2}$\\
\hline
\end{tabular}}
},
{(140,-10) *{\begin{tabular}{|c|}
\hline
$\frac32+\frac32$\\
\hline
$(1)_2 11, (2)_2$\\
$H_\mathrm{KSs}^{2+\frac43}$\\
\hline
\end{tabular}}
},
{\ar (150,-11);(166,-26)},
{(140,-38) *{\begin{tabular}{|c|}
\hline
$2+\frac43$\\
\hline
$(1)_3 1, (2)(2)$\\
$H_\mathrm{KSs}^{2+\frac43}$\\
\hline
\end{tabular}}
},
{\ar (150,-39);(166,-26)},
{\ar (150,-39);(165,-56)},
{(140,-65) *{\begin{tabular}{|c|}
\hline
$\frac54+1+1$\\
\hline
$(1)_4, 22, 22$\\
$H_\mathrm{KSs}^{2+\frac43}$\\
\hline
\end{tabular}}
},
{\ar (149,-66);(165,-56)},
{(175,-25) *{\begin{tabular}{|c|}
\hline
$\frac32+\frac43$\\
\hline
$(1)_3 1, (2)_2$\\
$H_\mathrm{KSs}^{2+\frac54}$\\
\hline
\end{tabular}}
},
{\ar (184,-26);(202,-38.5)},
{(175,-55) *{\begin{tabular}{|c|}
\hline
$2+\frac54$\\
\hline
$(1)_4, (2)(2)$\\
$H_\mathrm{KSs}^{2+\frac54}$\\
\hline
\end{tabular}}
},
{\ar (184,-56);(202,-38.5)},
{(210,-38) *{\begin{tabular}{|c|}
\hline
$\frac32+\frac54$\\
\hline
$(1)_4, (2)_2$\\
$H_\mathrm{KSs}^{\frac32+\frac54}$\\
\hline
\end{tabular}
}},
\end{xy}
\endgroup
}
\end{landscape}

\bigskip 

\noindent 
\textbf{Acknowledgements} 

\noindent 
The author wishes to thank Professors Hidetaka Sakai, Akane Nakamura, and
Kazuki Hiroe for their invaluable suggestions and comments.

\section{HTL canonical forms and their degeneration}\label{sec:HTL}
\subsection{HTL canonical forms}\label{sec:HTLform}
Here we recall HTL canonical forms of systems of linear differential equations.
Consider a linear system
\[
\frac{dY}{dx}=A(x)Y.
\]
The transformation $A(x) \mapsto A^P(x):=P(x)^{-1}A(x)P(x)-P(x)^{-1}\frac{dP(x)}{dx}$ by an invertible matrix $P(x)$ is called a gauge transformation.
This corresponds to the change of the dependent variable $Y=P(x)Z$.
Then $Z$ satisfies
\[
\frac{dZ}{dx}=A^P(x) Z.
\]

Linear differential equations that we treat in this series of papers are systems with rational function coefficients:
\begin{equation}
\frac{dY}{dx}=
\left(\sum_{\nu=1}^n\sum_{k=0}^{r_{\nu}}\frac{A_{\nu}^{(k)}}{(x-u_{\nu})^{k+1}}
+\sum_{k=1}^{r_{\infty}}A_{\infty}^{(k)}x^{k-1}
\right)Y, \quad
A_j^{(k)} \in M_m(\mathbb{C}).
\end{equation}
This system has singular points at $x=u_1,\ldots,u_n$, and $\infty$.
Taking the local coordinates $z=x-u_\nu$ or $z=1/x$,
we can write the system at $z=0$ as follows:
\begin{equation}\label{eq:Laurent}
\frac{dY}{dz}=\left(
\frac{A_0}{z^{r+1}}+\frac{A_1}{z^{r}}+\cdots+A_{r+1}+A_{r+2} z+\cdots
\right)Y.
\end{equation}
We denote the field of formal Laurent series in $z$ by $\mathbb{C}(\!(z)\!)$,
and the field of Puiseux series $\cup_{p > 0}\mathbb{C}(\!(z^{\frac1p})\!)$ by $\mathcal{K}_z$.

Here we give the definition of HTL canonical forms.
\begin{df}\label{df:HTL}
An element in $M_m(\mathcal{K}_z)$ of the form
\begin{equation}
\frac{D_0}{z^{l_0}}+\frac{D_1}{z^{l_1}}+\cdots+\frac{D_{s-1}}{z^{l_{s-1}}}+\frac{\Theta}{z^{l_s}}
\end{equation}
where
	\begin{itemize}
		\item $l_0, \ldots, l_{s}\ (l_0 > l_1 > \cdots > l_{s-1}>l_s = 1)$ are rational numbers,
		
		\item $D_0, \ldots, D_{s-1}$ are diagonal matrices,
		
		\item $\Theta$ is a (not necessarily diagonal) Jordan matrix which commutes with all $D_j$'s,
	\end{itemize}
is called an \textit{HTL canonical form} (or \textit{HTL form} for short).
\qed
\end{df}
Then the following theorem holds.
\begin{thm}[Hukuhara \cite{Huk}, Turrittin \cite{Tur}, Levelt \cite{Lev}]\label{thm:HTL}
For any element $A(z)$ in $M_m(\mathbb{C}(\!(z)\!))$:
\begin{equation}\label{eq:Laurent2}
A(z)=\frac{A_0}{z^{r+1}}+\frac{A_1}{z^{r}}+\cdots+A_{r+1}+A_{r+2} z+\cdots, \quad
A_j \in M_m(\mathbb{C}),
\end{equation}
there exists $P(z) \in \mathrm{GL}_m(\mathcal{K}_z)$ such that $A^P(z)$ is an HTL form:
\begin{equation}\label{eq:HTLform}
A^P(z)=
\frac{D_0}{z^{l_0}}+\frac{D_1}{z^{l_1}}+\cdots+\frac{D_{s-1}}{z^{l_{s-1}}}+\frac{\Theta}{z^{l_s}}.
\end{equation}
%
%
Here $l_0, \ldots, l_{s}$ are uniquely determined only by the original system (\ref{eq:Laurent2}).

If the following
\begin{equation}
\frac{\tilde{D}_0}{z^{l_0}}+\frac{\tilde{D}_1}{z^{l_1}}+\cdots+\frac{\tilde{D}_{s-1}}{z^{l_{s-1}}}
+\frac{\tilde{\Theta}}{z^{l_s}}
\end{equation}
is another HTL canonical form corresponding to the same system,
there exist a constant matrix $g \in \mathrm{GL}_m(\mathbb{C})$ and a natural number $k \in \mathbb{Z}_{\ge 1}$ such that
\begin{equation}
\tilde{D}_j=g^{-1}D_j g, \quad \exp(2\pi i k \tilde{\Theta})=g^{-1}\exp(2\pi i k \Theta)g
\end{equation}
hold.
\end{thm}
We call (\ref{eq:HTLform}) the HTL canonical form (or HTL form) of (\ref{eq:Laurent2}).
The number $l_0-1$ is called the {\it Poincar\'e rank} of the singular point.
When there is a rational number $l_j$ that is not an integer, the singular point is called a {\it ramified} irregular singular point.
A linear system is said to be of \textit{ramified type} if it has a ramified irregular singular point.


HTL forms can be computed by means of block diagonalizations and shearing transformations.
Here the block diagonalization means the following
\begin{prop}{(block diagonalization)}\label{thm:block_diag}
Let $A(z)$ be a formal Laurent series
\begin{equation}\label{eq:irreg_sys}
A(z)=\frac{A_0}{z^{r+1}}+\frac{A_1}{z^r}+\cdots
\quad (A_j \in M_m(\mathbb{C})).
\end{equation}
We write the eigenvalues of $A_0$ as $\lambda_1,\ldots,\lambda_n$ and their multiplicities $m_1,\ldots,m_n$ respectively.
Then we can choose a formal power series $P(z)$ so that the gauge transformation by $P(z)$ reduces {\rm (\ref{eq:irreg_sys})} to
the following form:
\begin{align}
A^P(z)
=
\begin{pmatrix}
B_1 & & \\
& \ddots & \\
& & B_n
\end{pmatrix}, \quad
B_k=\frac{B^k_0}{z^{r+1}}+\frac{B^k_1}{z^r}+\cdots
\end{align}
where $B^k_0=\lambda_k I_{m_k}+N_k$ with $N_k$ being nilpotent.
\end{prop}

We briefly explain shearing transformations in the next subsection.
In this series of papers we use the notations what we call
\textit{Riemann schemes} and \textit{spectral types} to express linear systems.
The Riemann scheme and the spectral type of a linear system
are defined by the HTL forms; see \cite{K1}. 

\subsection{Shearing transformations}
Here we illustrate how the shearing transformations work, following \cite{Wa}.
By virtue of Proposition~\ref{thm:block_diag},
it suffices to consider the following system
\begin{align*}
\frac{dY}{dz}&=A(z)Y, \\
A(z)&=\frac{1}{z^{r+1}}\left( A_0+A_1z+\cdots+A_r z^r+\cdots \right)
\quad (A_j \in M_m(\mathbb{C}))
\end{align*}
where $A_0$ is a Jordan matrix with only one eigenvalue.
By means of a scalar gauge transformation, we can shift $A_0$ by a scalar matrix.
Thus, without loss of generality, we can assume that $A_0$ is nilpotent.

Let $\tilde{A}(z)=A_0+A_1z+\cdots$ and write the $(i,j)$-entry of $\tilde{A}(z)$ as $a_{ij}(z)$.
We write a non-zero $a_{ij}(z)$ as
\begin{equation}
a_{ij}(z)=z^{\alpha_{ij}}a^*_{ij}(z), \quad
a^*_{ij}(0) \ne 0,
\end{equation}
where $\alpha_{ij}$ is a non-negative integer.
Note that at least one of the $\alpha_{i,i+1}$'s $(i=1,\ldots,m-1)$ is zero.

Now let $s$ be a positive constant and consider the gauge transformation by
\[
S=\mathrm{diag}(1, z^s, \ldots, z^{(m-1)s}).
\]
We denote the $(i,j)$-entry of 
$z^{r+1}A^S(z)$ by $b_{ij}(z)$.
If $a_{ij}(z) \neq 0$, then
\begin{equation}
b_{ij}(z)=z^{(j-i)s+\alpha_{ij}}a^*_{ij}(z)-\delta_{ij}(i-1)s z^r.
\end{equation}
We can write this (in a similar way as above) as $b_{ij}(z)=z^{\beta_{ij}}b^*_{ij}(z) \ (b^*_{ij}(0) \ne 0)$ again.
It is easy to see that $\beta_{ij}$ is a linear function in $s$.
Draw the graphs of $\beta=\beta_{ij}(s)$ in $(s, \beta)$-plane.
Consider the intersection points between the lines with negative slopes and the line $\beta=s$.
Find the intersection point with smallest $\beta$-coordinate 
and write it as $s_0$.
Then let $s=s_0$ in the shearing matrix $S$.
We repeat this procedure so that the coefficient matrix of the leading term becomes diagonalizable.
\begin{eg}
We take the system of $x$-direction of (\ref{eq:(1)_4,(2)_2}) as an example.
We compute the HTL form at $x=\infty$ of the system.
Changing the independent variable from $x$ to $z=1/x$ and modifying the gauge of the system so that $A_0$ is in Jordan canonical form,
we have the following system:
\begin{equation}
\frac{dY}{dz}=\frac{1}{z^2}\left( A_0+A_1z+A_2z^2 \right)Y,
\end{equation}
{\small
\begin{align*}
A_0&=
\begin{pmatrix}
0 & 1 & 0 & 0 \\
0 & 0 & 0 & 0 \\
0 & 0 & 0 & 0\\
0 & 0 & 0 & 0
\end{pmatrix}, \quad
A_1=
\begin{pmatrix}
-p_2q_2-\frac12 & -q_2/t & 0 & -\frac{q_1q_2}{t} \\
0 & p_2q_2+\frac12 & q_1 & 0 \\
0 & -q_2/t & -p_1q_1 & t/q_1 \\
1 & 0 & -1 & p_1q_1
\end{pmatrix}, \quad
A_2=
\begin{pmatrix}
0 & 0 & 0 & 0 \\
-t/q_2 & 0 & 0 & 0 \\
0 & 0 & 0 & -1 \\
0 & 0 & 0 & 0
\end{pmatrix}.
\end{align*}
}
The entries of $\tilde{A}(z)=A_0+A_1z+A_2z^2$ are explicitly written as follows:
\begin{center}
\begin{tabular}{llll}
$a_{11}=z(-p_2q_2-\frac12)$, & $a_{12}=1-\frac{q_2}{t}z$, & $a_{13}=0$, & $a_{14}=z(-\frac{q_1q_2}{t})$, \\
$a_{21}=z^2(-t/q_2)$, & $a_{22}=z(p_2q_2+\frac12)$, & $a_{23}=z(q_1)$, & $a_{24}=0$, \\
$a_{31}=0$, & $a_{32}=z( -q_2/t )$, & $a_{33}=z(-p_1q_1)$, & $a_{34}=z(t/q_1-z)$, \\
$a_{41}=z$, & $a_{42}=0$, & $a_{43}=-z$, & $a_{44}=z(p_1q_1)$.
\end{tabular}
\end{center}
Thus we have
\begin{center}
\begin{tabular}{llll}
$\beta_{11}=1$, & $\beta_{12}=s$, &  & $\beta_{14}=3s+1$, \\
$\beta_{21}=-s+2$, & $\beta_{22}=1$, & $\beta_{23}=s+1$, &  \\
                        & $\beta_{32}=-s+1$, & $\beta_{33}=1$, & $\beta_{34}=s+1$, \\
$\beta_{41}=-3s+1$, &  & $\beta_{43}=-s+1$, & $\beta_{44}=1$.
\end{tabular}
\end{center}
In this case, the $\beta$-coordinate of the intersection point of the lines $\beta=s$ and $\beta=-3s+1$ gives $s_0$.
Thus we find $s_0=1/4$.

Then we have
\begin{align}
A^S(z)=
{\small
\begin{pmatrix}
0 & 1 & 0 & 0 \\
0 & 0 & 0 & 0 \\
0 & 0 & 0 & 0 \\
1 & 0 & 0 & 0
\end{pmatrix}
\frac{1}{z^{\frac74}}+
\begin{pmatrix}
0 & 0 & 0 & 0 \\
0 & 0 & 0 & 0 \\
0 & -q_2/t & 0 & 0 \\
0 & 0 & -1 & 0
\end{pmatrix}
\frac{1}{z^{\frac54}}\cdots.}
\end{align}
Note that its leading coefficient has the following Jordan canonical form
{\small
\begin{equation}
\begin{pmatrix}
0 & 1 & 0 & 0 \\
0 & 0 & 0 & 0 \\
0 & 0 & 0 & 0 \\
1 & 0 & 0 & 0
\end{pmatrix}
\sim
\begin{pmatrix}
0 & 1 & 0 & 0 \\
0 & 0 & 1 & 0 \\
0 & 0 & 0 & 0 \\
0 & 0 & 0 & 0
\end{pmatrix}.
\end{equation}
}

Repeating such a procedure two more times, eventually we have
{\small
\begin{align*}
&\begin{pmatrix}
0 & 1 & 0 & 0 \\
0 & 0 & 1 & 0 \\
0 & 0 & 0 & 1 \\
t & 0 & 0 & 0
\end{pmatrix}
\frac{1}{z^{\frac54}}+\cdots
\\ \sim
&\begin{pmatrix}
t^{\frac14} & 0 & 0 & 0 \\
0 & \sqrt{-1}t^{\frac14} & 0 & 0 \\
0 & 0 & -t^{\frac14} & 0 \\
0 & 0 & 0 & -\sqrt{-1}t^{\frac14}
\end{pmatrix}
\frac{1}{z^{\frac54}}+
\begin{pmatrix}
-9/8 & * & * & * \\
* & -9/8 & * & * \\
* & * & -9/8 & * \\
* & * & * & -9/8
\end{pmatrix}
\frac{1}{z}+\cdots.
\end{align*}
}
Here the leading coefficient becomes diagonalizable.
Then we obtain the following HTL form at $x=\infty$ of the system of $x$-direction of (\ref{eq:(1)_4,(2)_2}):
{\small
\begin{align*}
\begin{pmatrix}
t^{\frac14} & 0 & 0 & 0 \\
0 & \sqrt{-1}t^{\frac14} & 0 & 0 \\
0 & 0 & -t^{\frac14} & 0 \\
0 & 0 & 0 & -\sqrt{-1}t^{\frac14}
\end{pmatrix}
\frac{1}{z^{\frac54}}+
\begin{pmatrix}
-9/8 & 0 & 0 & 0 \\
0 & -9/8 & 0 & 0 \\
0 & 0 & -9/8 & 0 \\
0 & 0 & 0 & -9/8
\end{pmatrix}
\frac{1}{z}.
\end{align*}
}
This can be seen as the direct sum of $\frac{t^{1/4}}{z^{5/4}}+\frac{-9/8}{z}$ and its copies.
Here by ``copies'' we mean the images of $\frac{t^{1/4}}{z^{5/4}}+\frac{-9/8}{z}$ under the cyclic action $z^{1/4} \mapsto \sqrt{-1}z^{1/4}$.
Moreover we can eliminate $-\frac{9/8}{z} I$ by means of the gauge transformation by $z^{-\frac98}$.
Thus we can express the HTL form as
\[
\begin{array}{c}
  x=\infty \, \left( \frac14 \right) \\
\overbrace{
\begin{array}{cc}
t^{\frac14} & 0 \\
\sqrt{-1}t^{\frac14}  & 0 \\
-t^{\frac14}  & 0 \\
-\sqrt{-1}t^{\frac14} & 0
\end{array}}
\end{array}.
\]
As explained in \cite{K1},
we write the spectral type of the singular point $x=\infty$ of this system as $(1)_4$.
\qed
\end{eg}

\subsection{Degeneration of HTL forms}\label{sec:deg_HTL}
In the case of  the original Painlev\'e equations (with standard $2 \times 2$ Lax pairs),
the degenerations of HTL forms are caused by the degenerations of Jordan canonical forms
of the leading coefficient matrix at irregular singular points.
However, in general, degeneration of an HTL form does not always correspond to degeneration of a Jordan canonical form.

In this subsection, we make some comments on degeneration of HTL forms.
For detailed calculations, see Appendix~\ref{sec:appendix}.

We begin with a simple example.
Let $A(x)$ be the following Laurent series:
\begin{equation}
A(x)=\frac{A_0}{x^2}+\frac{A_1}{x}+\cdots
\end{equation}
where
\[
A_0=
\begin{pmatrix}
0 & 1 \\
0 & 0
\end{pmatrix}, \quad
A_k=\left( a^{(k)}_{ij} \right) \in M_2(\mathbb{C}).
\]
Let $S=\mathrm{diag}(1, x^{1/2})$. Then we have
\begin{equation}
A^S(x)=
\frac{1}{x^{\frac32}}
\begin{pmatrix}
0 & 1 \\
a^{(1)}_{21} & 0
\end{pmatrix}+
\frac{1}{x}
\begin{pmatrix}
a^{(1)}_{11} & 0 \\
0 & a^{(1)}_{22}-\frac12
\end{pmatrix}+\cdots.
\end{equation}
Here the coefficient matrix of the leading term becomes diagonalizable
provided that $a^{(1)}_{21} \neq 0$,
and we find that the corresponding system of linear differential equations
has $x=0$ as an irregular singular point of
Poincar\'e rank 1/2.
The HTL form at $x=0$ is
\[
\begin{array}{c}
  x=0 \, \left( \frac12 \right) \\
\overbrace{
\begin{array}{cc}
\sqrt{a^{(1)}_{21}} & \mathrm{tr}(A_1) \\
-\sqrt{a^{(1)}_{21}}  & \mathrm{tr}(A_1)
\end{array}}
\end{array}.
\]
If $a^{(1)}_{21}=0$, using $\tilde{S}=\mathrm{diag}(1, x)$ instead of the above $S$,
we have the different HTL form (we omit the details).
This example shows us that if the leading matrix is a Jordan canonical form whose $(1,2)$-entry is one,
then whether or not the $(2,1)$-entry of the next matrix is zero is meaningful.

Now we consider the degeneration $(111)(1) \to (1)_2 11$.
We take the degeneration $(111)(1), (2)(2) \to (1)_2 11, (2)(2)$ as an example.
Note that the linear system $(111)(1), (2)(2)$ is of the form~\cite{KNS}:
\begin{equation*}
\frac{dY}{dx}=\left(
\frac{A_0^{(-1)}}{x^2}+\frac{A_0^{(0)}}{x}+A_\infty
\right)Y
\end{equation*}
where
\[
A_\infty=\mathrm{diag}(-t,0,0,0).
\]
According to Appendix~\ref{sec:appendix}, there exists an invertible matrix, say $G$, such that
\begin{equation*}
G^{-1}A_\infty G=
{\small
\begin{pmatrix}
0 & 1 & 0 & 0 \\
0 & 0 & 0 & 0 \\
0 & 0 & 0 & 0 \\
0 & 0 & 0 & 0
\end{pmatrix}}
+O(\varepsilon) \quad
(\varepsilon \to 0).
\end{equation*}
In other words, the semisimple matrix $A_\infty$ degenerates to a nilpotent matrix as $\varepsilon$ tends to zero.
The degeneration of this HTL form corresponds to the degeneration of the Jordan canonical form.

However, in the case of the degeneration $(1)_2 11 \to (1)_3 1$ (see, for example, $(1)_211, (2)(2) \to (1)_31, (2)(2)$ in Appendix~\ref{sec:appendix}),
the Jordan canonical form of the leading coefficient does not change.
Instead, the $(2,1)$-entry of the next matrix becomes zero,
so that the HTL form $(1)_2 11$ degenerates to $(1)_3 1$.


\section{Lax pairs of degenerate Sasano systems}\label{sec:degeneration}
%

The Sasano system was discovered by Sasano~\cite{Ss} in his study on a generalization of Painlev\'e equations
from the viewpoint of spaces of initial values,
and later it was obtained as a similarity reduction of the Drinfel'd-Sokolov hierarchy \cite{FS2}.

Sakai \cite{Sak2} derived the Sasano system from the isomonodromic deformation of the following Fuchsian system:
\begin{equation}\label{eq:Fuchs_sasano}
\frac{dY}{dx}=
\left(
\frac{A_0}{x}+\frac{A_1}{x-1}+\frac{A_t}{x-t}
\right)Y
\end{equation}
where $A_0$, $A_1$, and $A_t$ are $4 \times 4$ matrices satisfying the following conditions
\begin{equation}\label{eq:eigen_condition}
A_0 \sim \mathrm{diag}(0,0,0,\theta^0), \ 
A_1 \sim \mathrm{diag}(0,0,\theta^1,\theta^1), \ 
A_t \sim \mathrm{diag}(0,0,\theta^t,\theta^t), \ 
\end{equation}
and
\begin{equation}\label{residue_infty}
A_\infty:=-(A_0+A_1+A_t)=
\mathrm{diag}(\theta^\infty_1,\theta^\infty_2,\theta^\infty_3,\theta^\infty_4).
\end{equation}
Thus the spectral type of the Fuchsian system (\ref{eq:Fuchs_sasano}) is $31,22,22,1111$.
Taking the trace of (\ref{residue_infty}), we have the Fuchs relation
\[
\theta^0+2(\theta^1+\theta^t)+\theta^\infty_1+\theta^\infty_2+\theta^\infty_3+\theta^\infty_4=0.
\]

The isomonodromic deformation equation of (\ref{eq:Fuchs_sasano}) is equivalent to the Hamiltonian system
\begin{equation}\label{eq:Ss_system}
\frac{dq_i}{dt}=\frac{\partial H_{\mathrm{Ss}}^{D_6}}{\partial p_i}, \quad
\frac{dp_i}{dt}=-\frac{\partial H_{\mathrm{Ss}}^{D_6}}{\partial q_i} \quad (i=1,2)
\end{equation}
where the Hamiltonian is given by
\begin{align}\label{eq:Ss_Ham}
H_{\mathrm{Ss}}^{D_6}\left({\alpha, \beta, \gamma \atop \delta, \epsilon, \zeta};t;
{q_1,p_1 \atop q_2,p_2}\right)&=
H_{\mathrm{VI}}\left({\beta+\gamma+2\delta+\epsilon+\zeta, -\beta-\zeta \atop -\beta-2\gamma-2\delta-\epsilon,1-\alpha-\beta-2\delta-\epsilon-\zeta};t;q_1,p_1\right)\\
&\quad+H_{\mathrm{VI}}\left({\gamma+\delta, \epsilon \atop \zeta, 1-\alpha-\gamma};t;q_2,p_2\right)\nonumber\\
&\quad+\frac{2}{t(t-1)}(q_1-1)p_{2} q_{2}\{(q_{1}-t)p_{1}-\beta-\gamma-2\delta-\epsilon-\zeta)\} .\nonumber
\end{align}
We call the Hamiltonian system (\ref{eq:Ss_system}), (\ref{eq:Ss_Ham}) the Sasano system.
The parametrization of the Fuchsian system (\ref{eq:Fuchs_sasano}) is rather complicated, see \cite{KNS} for details.


In the present paper we obtain the following Hamiltoians, which did not appear in \cite{KNS}.
To our knowledge, these systems have not appeared in the literature:
{\allowdisplaybreaks
\begin{align}
tH_\mathrm{KSs}^{2+\frac32}\left( {\alpha, \beta \atop \gamma}; t; {q_1,p_1\atop q_2,p_2} \right)&=
tH_{\mathrm{III}(D_7)}\left( \alpha; t ; q_1, p_1 \right)
+tH_{\mathrm{III}(D_7)}\left( \beta; t; q_2, p_2 \right)
+2p_2q_1(p_1q_1+\gamma)-q_1, \\
tH_\mathrm{KSs}^{2+\frac43}\left( {\alpha, \beta}; t; {q_1, p_1 \atop q_2, p_2} \right)&=
tH_{\mathrm{III}(D_7)}\left( \alpha; t; q_1, p_1 \right)
+tH_{\mathrm{III}(D_7)}\left( \beta; t; q_2, p_2 \right)
-t(2p_1p_2+p_1+p_2), \\
tH_\mathrm{KSs}^{2+\frac54}\left( \alpha; t; {q_1, p_1 \atop q_2, p_2} \right)&=
tH_{\mathrm{III}(D_7)}\left( \alpha; t ;q_1, p_1 \right)
+tH_{\mathrm{III}(D_7)}\left( \alpha; t; q_2, p_2 \right)
-t\left( 2\frac{p_2}{q_1}+p_1+p_2 \right), \\
tH_\mathrm{KSs}^{\frac32+\frac54}\left( t; {q_1, p_1 \atop q_2, p_2} \right)&=
tH_{\mathrm{III}(D_8)}\left(t ;q_1,p_1\right)+tH_{\mathrm{III}(D_8)}\left(t; q_2,p_2\right)
-2\frac{q_1q_2}{t}+q_1+q_2.
\end{align}
}
\begin{rem}
Although the Sasano system is a system of ordinary differential equations,
there are four systems of partially differential equations in the degeneration scheme in Section~\ref{sec:intro}:
$H_{\mathrm{Gar}}^{\frac32+1+1+1}$ (associated with $(1)_2(1)_2, 31, 22$), $H_{\mathrm{Gar}}^{\frac52+1+1}$ (associated with $(((1)(1)))_2, 31$),
$H_{\mathrm{Gar}}^{2+\frac32+1}$ (associated with $((1))_4, 31, 22$), and $H_{\mathrm{Gar}}^{\frac52+2}$ (associated with $((((((1))))))_4, 31$).
This should be interpreted as follows.

For example, the linear system $(11)(11), 31, 22$ degenerates to $(1)_2(1)_2, 31, 22$.
The former admit one-dimensional deformation, while the latter naturally has two-dimensional deformation.
By considering full deformation of the $(1)_2(1)_2, 31, 22$ system, we have the Hamltonians for
the degenerate Garnier system of type $3/2+1+1+1$.
The same argument applies to the other systems.
Therefore there are degenerate Garnier systems in two variables in the degeneration scheme.
\qed
\end{rem}

\subsection{Singularity pattern $\frac32+1+1$}
\subsubsection{Spectral type $(2)_2, 31, 1111$}
The Riemann scheme is given by
\[
\left(
\begin{array}{ccc}
  x=0 & x=1 \, \left( \frac12 \right) & x=\infty \\
\begin{array}{c} 0 \\ 0 \\ 0 \\ \theta^0 \end{array} &
\overbrace{\begin{array}{cc}
\sqrt{t} & 0\\
\sqrt{t} & 0\\
-\sqrt{t} & 0\\
-\sqrt{t} & 0
      \end{array}}&
\begin{array}{c} \theta^\infty_1 \\ \theta^\infty_2 \\ \theta^\infty_3\\ \theta^\infty_4 \end{array}
\end{array}
\right) ,
\]
and the Fuchs-Hukuhara relation is written as
$\theta^0+\theta^\infty_1+\theta^\infty_2+\theta^\infty_3+\theta^\infty_4=0$.

The Lax pair is expressed as
\begin{equation}\label{eq:Lax(2)_2,31,1111}
\left\{
\begin{aligned}
\frac{\partial Y}{\partial x}&=
\left(
\frac{A_0}{(x-1)^2}+\frac{A_1}{x-1}+\frac{A_2}{x}
\right)Y, \\
\frac{\partial Y}{\partial t}&=-\frac{1}{x-1}\left(\frac{A_0}{t}\right)Y,
\end{aligned}
\right.
\end{equation}
where
{\allowdisplaybreaks
\begin{align*}
A_{\xi}&=
U^{-1}
P^{-1}\hat{A}_\xi PU, \quad
P=
\begin{pmatrix}
  1 & 0 & 0 & 0 \\
 \frac{a_{21}}{\theta^\infty_2-\theta^\infty_1} & 1 & 0 & 0 \\
 \frac{a_{31}}{\theta^\infty_3-\theta^\infty_1} & 0 & 1 & 0 \\
 \frac{a_{41}}{\theta^\infty_4-\theta^\infty_1} & 0 & 0 & 1 
\end{pmatrix},\quad
U=\mathrm{diag}(1,u,v,w), \\
\hat{A}_0&=
\begin{pmatrix}
I_2 \\
\hat{B}
\end{pmatrix}
\begin{pmatrix}
\hat{C}\hat{B} &-\hat{C}
\end{pmatrix},\quad
\hat{B}=
\begin{pmatrix}
-q_1 & f_2\\
f_1q_1-q_2 & f_3
\end{pmatrix},\quad
\hat{C}=
\begin{pmatrix}
f_1 & 1\\
f_1-1 & 1
\end{pmatrix},
\\
\hat{A}_1&=
\begin{pmatrix}
-\theta^0-\theta^\infty_1 & a_{12} & (1-p_2)f_1-p_1 & 1-p_2\\
a_{21} & -\theta^\infty_2 & 0 & 0\\
a_{31} & 0 & -\theta^\infty_3& 0 \\
a_{41} & 0 & 0 & -\theta^\infty_4
\end{pmatrix},\\
\hat{A}_2&=
\begin{pmatrix}
1\\
0\\
0\\
0
\end{pmatrix}
\begin{pmatrix}
\theta^0 & -a_{12} & (p_2-1)f_1+p_1 & p_2-1
\end{pmatrix},
\end{align*}
}
and
{\allowdisplaybreaks
\begin{align*}
a_{12}&= (p_2-1)(f_1f_2+f_3)+p_1(f_2-q_1)-\theta^\infty_2-\theta^\infty_3,\quad
a_{21}=(p_2-1)(q_2-q_1)+\theta^\infty_2+\theta^\infty_4, \\
a_{31}&=q_1(-p_1q_1-p_2q_2+q_2+\theta^0+\theta^\infty_1-\theta^\infty_3)+((p_2-1)(q_2-q_1)+\theta^\infty_2+\theta^\infty_4)f_2-t,\\
a_{41}&=((\theta^\infty_3-\theta^\infty_4)f_3+t)(f_1-1)+(p_1+p_2-1)(q_2-q_1)f_3\\
&\quad+(f_1q_1-q_2)(p_1q_1+(p_2-1)q_2-\theta^0-\theta^\infty_1+\theta^\infty_4), \\
f_1&=\frac{1}{\theta_3^\infty-\theta_4^\infty}\{ p_1(q_1-q_2)+\theta^\infty_2+\theta^\infty_3 \},\quad
f_2=\frac{1}{\theta^\infty_2-\theta^\infty_3}\{ (p_1q_1+\theta^\infty_2+\theta^\infty_3)q_1+t \},\\
f_3&=\frac{1}{\theta^\infty_2-\theta^\infty_4}\{ (q_2-q_1f_1)(p_1q_1+\theta^\infty_2+\theta^\infty_3)-t f_1 \}.
\end{align*}
}

The Hamiltonian is given by
\begin{align}
&tH_{\mathrm{Ss}}^{D_4}
\left({-2\theta^\infty_1+1,\theta^\infty_1+\theta^\infty_2+\theta^\infty_3+\theta^\infty_4
\atop \theta^\infty_1-\theta^\infty_4-1, \theta^\infty_1-\theta^\infty_3-1};t;
{q_1,p_1\atop q_2,p_2}\right)\\
\quad&=tH_{\mathrm{III}(D_6)}
\left({\theta^\infty_2+\theta^\infty_3, 2\theta^\infty_3+1};t;q_1,p_1
\right)+tH_{\mathrm{III}(D_6)}\left({\theta^\infty_4-\theta^\infty_1+1, 2\theta^\infty_4+1};t;q_2,p_2\right)\nonumber\\
&\quad+2p_2q_1(p_1q_1+\theta^\infty_2+\theta^\infty_3).\nonumber
\end{align}

The gauge parameters $u$, $v$, $w$ satisfy
\begin{align}
\frac{t}{u}\frac{du}{dt}&=(2p_2-1)q_1-q_2,\quad
\frac{t}{v}\frac{dv}{dt}=2q_1(p_1+p_2)-q_1-q_2+2\theta^\infty_3+1,\\
\frac{t}{w}\frac{dw}{dt}&=2q_2(p_2-1)+2\theta^\infty_4+1.\nonumber
\end{align}

\subsubsection{Spectral type $(1)_2 11, 22, 22$}
The Riemann scheme is given by
\[
\left(
\begin{array}{ccc}
  x=0 & x=1 & x=\infty \, \left( \frac12 \right) \\
\begin{array}{c} 0 \\ 0 \\ \theta^0 \\ \theta^0 \end{array} &
\begin{array}{c} 0 \\ 0 \\ \theta^1 \\ \theta^1 \end{array} &
\overbrace{
\begin{array}{cc}
0 & \theta^\infty_2 \\
0 &  \theta^\infty_3 \\
\sqrt{t} & \theta^\infty_1/2\\
-\sqrt{t} & \theta^\infty_1/2
\end{array}}
\end{array}
\right) ,
\]
and the Fuchs-Hukuhara relation is written as
$2\theta^0+2\theta^1+\theta^\infty_1+\theta^\infty_2+\theta^\infty_3=0$.

The Lax pair is expressed as
\begin{equation}
\left\{
\begin{aligned}
\frac{\partial Y}{\partial x}&=
\left(
\frac{A_0}{x}+\frac{A_1}{x-1}+A_\infty
\right)Y, \\
\frac{\partial Y}{\partial t}&=\left( \frac1t A_\infty x+B \right)Y.
\end{aligned}
\right.
\end{equation}
Here
{\allowdisplaybreaks
\begin{align*}
A_\infty&=
\begin{pmatrix}
0 & t & 0 & 0 \\
0 & 0 & 0 & 0 \\
0 & 0 & 0 & 0 \\
0 & 0 & 0 & 0
\end{pmatrix}, \quad
A_0=
\begin{pmatrix}
0 & 0 & 0 & 0 \\
1-p_1-p_2 & \theta^0 & 0 & p_1 \\
-1 & 0 & \theta^0 & 1 \\
0 & 0 & 0 & 0
\end{pmatrix}, \quad
A_1=
\begin{pmatrix}
\theta^1 I_2-C_1 B_1 \\
B_1
\end{pmatrix}
\begin{pmatrix}
I_2 & C_1
\end{pmatrix}, \\
B_1&=
\begin{pmatrix}
1 & -q_1 \\
p_2(q_2-q_1)+\theta^1+\theta^\infty_3 & 
p_2q_2(q_1-q_2)-\theta^\infty_3 q_1-\theta^1 q_2
\end{pmatrix},\\
C_1&=
\frac{1}{q_1-q_2}
\begin{pmatrix}
p_1q_1(q_2-q_1)+(\theta^0+\theta^1+\theta^\infty_1)q_1+(\theta^0+\theta^\infty_2)q_2 & q_2 \\
p_1(q_2-q_1)-\theta^1-\theta^\infty_3 & 1
\end{pmatrix},\\
B&=\frac{1}{t}
\begin{pmatrix}
p_1q_1+p_2q_2-\theta^0-\theta^\infty_1 & 0 & 0 & \theta^0+\theta^1+\theta^\infty_1-p_1q_1 \\
1 & -p_1q_1-p_2q_2+\theta^0 & B_{23} & 0 \\
0 & -q_1 & B_{33} & 0 \\
p_2(q_2-q_1)+\theta^1+\theta^\infty_3 & 0 & 0 & (2p_2-1)q_1-\theta^\infty_3 \\
\end{pmatrix},
\end{align*}
where
\begin{align*}
B_{23}&=p_1(p_2q_2-\theta^0-\theta^\infty_2)-p_2(p_1q_1-\theta^0-\theta^1-\theta^\infty_1),\\
B_{33}&=(2p_1+2p_2-1)q_1+\theta^0+\theta^1+\theta^\infty_2.
\end{align*}
}
The Hamiltonian is given by
\begin{align}
&tH_{\mathrm{Ss}}^{D_4}
\left({\theta^0+\theta^1, \theta^1+\theta^\infty_2+\theta^\infty_3 \atop
-\theta^1, -\theta^0-\theta^1-\theta^\infty_2}; t; {q_1, p_1 \atop q_2, p_2}\right)\\
&=tH_{\mathrm{III}(D_6)}\left(-\theta^0-\theta^1-\theta^\infty_1, \theta^0+\theta^1+2\theta^\infty_2; t; q_1, p_1\right)
+tH_{\mathrm{III}(D_6)}\left(\theta^1, \theta^1-\theta^0; t; q_2, p_2\right)\nonumber\\
&\quad+2p_2q_1(p_1q_1-\theta^0-\theta^1-\theta^\infty_1).\nonumber
\end{align}

\subsubsection{Spectral type $(1)_2(1)_2, 31, 22$}
The Riemann scheme is given by
\[
\left(
\begin{array}{ccc}
  x=0 & x=1 & x=\infty \, \left( \frac12 \right) \\
\begin{array}{c} 0 \\ 0 \\ 0 \\ \theta^{0} \end{array} &
\begin{array}{c} 0 \\ 0 \\ \theta^{1} \\ \theta^{1} \end{array} &
\overbrace{\begin{array}{cc}
\sqrt{-t_1} & \theta^\infty_1/2\\
-\sqrt{-t_1} & \theta^\infty_1/2\\
\sqrt{-t_2} & \theta^\infty_2/2\\
-\sqrt{-t_2} & \theta^\infty_2/2
      \end{array}}
\end{array}
\right) ,
\]
and the Fuchs-Hukuhara relation is written as
$\theta^0+2\theta^1+\theta^\infty_1+\theta^\infty_2=0$.

The Lax pair is expressed as
\begin{equation}
\left\{
\begin{aligned}
\frac{\partial Y}{\partial x}&=
\left(
\frac{A_0}{x}+\frac{A_1}{x-1}+A_{\infty}
\right)Y, \\
\frac{\partial Y}{\partial t_i}&=(N_i x+B_i)Y \quad (i=1,2).
\end{aligned}
\right.
\end{equation}
Here
{\allowdisplaybreaks
\begin{align*}
A_0&=P
\begin{pmatrix}
1 \\
p_1 \\
p_2
\end{pmatrix}
\begin{pmatrix}
1-p_1-p_2 & 1 & 1
\end{pmatrix}Q, \quad
A_1=P
\begin{pmatrix}
1 & 1 \\
p_1-1 & p_1\\
p_2 & p_2-1
\end{pmatrix}
\begin{pmatrix}
p_1 & -1 & 0 \\
p_2 & 0 & -1 
\end{pmatrix}Q,\\
A_\infty&=
\begin{pmatrix}
0 & 0 & t_1 & 0 \\
0 & 0 & 0 & t_2 \\
0 & 0 & 0 & 0 \\
0 & 0 & 0 & 0
\end{pmatrix}, \quad
N_1=
\begin{pmatrix}
0 & 0 & 1 & 0 \\
0 & 0 & 0 & 0 \\
0 & 0 & 0 & 0 \\
0 & 0 & 0 & 0
\end{pmatrix}, \quad
N_2=
\begin{pmatrix}
0 & 0 & 0 & 0 \\
0 & 0 & 0 & 1 \\
0 & 0 & 0 & 0 \\
0 & 0 & 0 & 0
\end{pmatrix}, \\
P&=
\begin{pmatrix}
(1-p_1-p_2)(p_1q_1-\theta^\infty_1)+\theta^1(p_2-1) & 0 & 0 \\
(1-p_1-p_2)(p_2q_2-\theta^\infty_2)+\theta^1(p_1-1) & 0 & 0 \\
0 & -1 & 0 \\
0 & 0 & -1
\end{pmatrix}, \\
Q&=
\begin{pmatrix}
0 & 0 & q_1 & q_2 \\
1 & 0 & \theta^\infty_1 & p_1(q_2-q_1)+\theta^1+\theta^\infty_1 \\
0 & 1 & p_2(q_1-q_2)+\theta^1+\theta^\infty_2 & \theta^\infty_2
\end{pmatrix},\\
B_1&=
\begin{pmatrix}
(B_1)_{11} & \frac{p_1(q_1-q_2)-\theta^1-\theta^\infty_1}{t_1-t_2} & (B_1)_{13} & p_1 \\
\frac{t_2(p_2(q_2-q_1)-\theta^1-\theta^\infty_2)}{t_1(t_1-t_2)} & (B_1)_{22} & (B_1)_{23} & 0 \\
-1/t_1 & 0 & (B_1)_{11}-\frac{\theta^\infty_1}{t_1} & \frac{p_1(q_1-q_2)-\theta^1-\theta^\infty_1}{t_1-t_2} \\
0 & 0 & \frac{p_2(q_2-q_1)-\theta^1-\theta^\infty_2}{t_1-t_2} & (B_1)_{44}
\end{pmatrix},\\
B_2&=
\begin{pmatrix}
(B_2)_{11} & \frac{t_1(p_1(q_1-q_2)-\theta^1-\theta^\infty_1)}{t_2(t_2-t_1)} & 0 & (B_2)_{14} \\
\frac{p_2(q_2-q_1)-\theta^1-\theta^\infty_2}{t_2-t_1} & (B_2)_{22} & p_2 & (B_2)_{24} \\
0 & 0 & (B_2)_{33} & \frac{p_1(q_1-q_2)-\theta^1-\theta^\infty_1}{t_2-t_1} \\
0 & -1/t_2 & \frac{p_2(q_2-q_1)-\theta^1-\theta^\infty_2}{t_2-t_1} & (B_2)_{22}-\frac{\theta^\infty_2}{t_2}
\end{pmatrix},
\end{align*}
}
where
{\allowdisplaybreaks
\begin{align*}
(B_1)_{11}&=\frac{p_1(q_1-q_2)-\theta^\infty_1}{t_1-t_2},\\
(B_1)_{13}&=\frac{q_1}{t_1}((p_1q_1-\theta^\infty_1)(1-p_1-p_2)+\theta^1(p_2-1))+p_1-1,\\
(B_1)_{22}&=\frac{1}{t_1}((p_1+p_2-1)q_1-\theta^\infty_1)+\frac{1}{t_1-t_2}(p_2(q_2-q_1)-\theta^\infty_2),\\
(B_1)_{23}&=\frac{q_1}{t_1}((p_2q_2-\theta^\infty_2)(1-p_1-p_2)+\theta^1(p_1-1)),\\
(B_1)_{44}&=(B_1)_{22},\\
(B_2)_{11}&=\frac{1}{t_2}((p_1+p_2-1)q_2-\theta^\infty_2)+\frac{1}{t_2-t_1}(p_1(q_1-q_2)-\theta^\infty_1),\\
(B_2)_{14}&=\frac{q_2}{t_2}((p_1q_1-\theta^\infty_1)(1-p_1-p_2)+\theta^1(p_2-1)),\\
(B_2)_{22}&=\frac{p_2(q_2-q_1)-\theta^\infty_2}{t_2-t_1},\\
(B_2)_{24}&=\frac{q_2}{t_2}((p_2q_2-\theta^\infty_2)(1-p_1-p_2)+\theta^1(p_1-1))+p_2-1,\\
(B_2)_{33}&=(B_2)_{11}.
\end{align*}
}

The Hamiltonians are given by
\begin{align}
&t_1H_{\mathrm{Gar}, t_1}^{\frac32+1+1+1}
 \left({\theta^1+\theta^\infty_1,\theta^1+\theta^\infty_2 \atop \theta^1};{t_1 \atop t_2};{q_1,p_1\atop q_2,p_2}\right)  \\
&\qquad=t_1H_\mathrm{III(D_6)}\left({-\theta^1-\theta^\infty_1, -\theta^\infty_1};t_1;q_1,p_1\right)+p_2q_1(p_1q_1-\theta^1-\theta^\infty_1)\nonumber\\
&\qquad \qquad+\frac{t_1}{t_1-t_2}(p_1(q_1-q_2)-\theta^1-\theta^\infty_1)(p_2(q_2-q_1)-\theta^1-\theta^\infty_2), \nonumber \\
&t_2H_{\mathrm{Gar}, t_2}^{\frac32+1+1+1}
 \left({\theta^1+\theta^\infty_1,\theta^1+\theta^\infty_2 \atop \theta^1};{t_1 \atop t_2};{q_1,p_1\atop q_2,p_2}\right)  \\
&\qquad=t_2H_\mathrm{III(D_6)}\left({-\theta^1-\theta^\infty_2, -\theta^\infty_2};t_2;q_2,p_2\right)+p_1q_2(p_2q_2-\theta^1-\theta^\infty_2)\nonumber\\
&\qquad \qquad+\frac{t_2}{t_2-t_1}(p_1(q_1-q_2)-\theta^1-\theta^\infty_1)(p_2(q_2-q_1)-\theta^1-\theta^\infty_2)\nonumber.
\end{align}

\subsubsection{Spectral type $((1))_4, 31, 22$}
The Riemann scheme is given by
\[
\left(
\begin{array}{ccc}
  x=0 & x=1 & x=\infty \, \left( \frac14 \right) \\
\begin{array}{c} 0 \\ 0 \\ 0 \\ \theta^{0} \end{array} &
\begin{array}{c} 0 \\ 0 \\ \theta^{1} \\ \theta^{1} \end{array} &
\overbrace{\begin{array}{ccc}
\sqrt{t_1} & \frac{1}{2}\sqrt{t_2}{t_1}^{-1/4} & \theta^\infty_1/4\\
-\sqrt{t_1} & \frac{\sqrt{-1}}{2}\sqrt{t_2}{t_1}^{-1/4} & \theta^\infty_1/4\\
\sqrt{t_1} & -\frac{1}{2}\sqrt{t_2}{t_1}^{-1/4} & \theta^\infty_1/4\\
-\sqrt{t_1} & -\frac{\sqrt{-1}}{2}\sqrt{t_2}{t_1}^{-1/4} & \theta^\infty_1/4
      \end{array}}
\end{array}
\right) ,
\]
and the Fuchs-Hukuhara relation is written as
$\theta^0+2\theta^1+\theta^\infty_1=0$.

The Lax pair is expressed as
\begin{equation}
\left\{
\begin{aligned}
\frac{\partial Y}{\partial x}&=
\left(
\frac{A_0}{x}+\frac{A_1}{x-1}+A_{\infty}
\right)Y, \\
\frac{\partial Y}{\partial t_i}&=(N_i x+B_i)Y \quad (i=1,2).
\end{aligned}
\right.
\end{equation}
Here
{\allowdisplaybreaks
\begin{align*}
A_0&=
\begin{pmatrix}
\theta^0 \\
-1 \\
1-p_1 \\
-q_2/t_2
\end{pmatrix}
\begin{pmatrix}
1 & 0 & 0 & 0
\end{pmatrix}, \quad
A_1=
\begin{pmatrix}
-t_2p_2 & -q_1(p_1q_1+\theta^1) \\
p_1q_1+\theta^1 & \frac{q_1q_2}{t_2}(p_1q_1+\theta^1)\\
0 & -p_1q_1 \\
1 & 0
\end{pmatrix}
\begin{pmatrix}
q_2/t_2 & 1 & 0 & -p_1q_1+p_2q_2 \\
-1/q_1 & 0 & 1 & -t_2p_2/q_1
\end{pmatrix},\\
A_\infty&=
\begin{pmatrix}
0 & 0 & t_1 & t_2 \\
0 & 0 & 0 & t_1 \\
0 & 0 & 0 & 0 \\
0 & 0 & 0 & 0
\end{pmatrix}, \quad
N_1=
\begin{pmatrix}
0 & 0 & 1 & 0 \\
0 & 0 & 0 & 1 \\
0 & 0 & 0 & 0 \\
0 & 0 & 0 & 0
\end{pmatrix}, \quad
N_2=
\begin{pmatrix}
0 & 0 & 0 & 1 \\
0 & 0 & 0 & 0 \\
0 & 0 & 0 & 0 \\
0 & 0 & 0 & 0
\end{pmatrix}, \\
B_1&=\frac{1}{t_1}
\begin{pmatrix}
0 & 0 & 0 & 0 \\
0 & -p_2q_2 & \frac{q_2}{t_2}(q_1(p_1q_1+\theta^1)-t_1-t_2p_2) & (1-p_1)(q_1(p_1q_1+\theta^1)-t_1)+t_2p_1p_2 \\
1 & -t_2/t_1 & -2p_1q_1+p_2q_2-\frac{t_2}{t_1}-\theta^0-\theta^1 & t_2p_2(p_1+1)+\frac{t_2}{t_1}(2p_1q_1+\theta^1) \\
0 & 1 & 1 & -2p_1q_1-\theta^1
\end{pmatrix},\\
B_2&=\frac{1}{t_2}
\begin{pmatrix}
0 & 0 & 0 & 0 \\
-1 & 2p_2q_2-\theta^0 & p_2q_2 & -t_2p_2(p_1+1) \\
0 & t_2/t_1 & t_2/t_1 & -\frac{t_2}{t_1}(2p_1q_1+\theta^1) \\
0 & 0 & -1 & 2p_2q_2-\theta^0
\end{pmatrix}.
\end{align*}
}

The Hamiltonians are given by
\begin{align}
&t_1H_{\mathrm{Gar},t_1}^{2+\frac{3}{2}+1}\left({-\theta^1, -\theta^0};{t_1 \atop t_2};{q_1,p_1 \atop q_2,p_2}\right)\\
&=
t_1H_{\mathrm{III}(D_6)}\left({\theta^1, \theta^0+\theta^1}; t_1; q_1,p_1\right)
+t_2p_1p_2-p_2q_2(2p_1q_1+\theta^1)+\frac{t_2}{t_1}p_1q_1-\frac{q_1q_2}{t_1}(p_1q_1+\theta^1),\nonumber \\
&t_2H_{\mathrm{Gar},t_2}^{2+\frac{3}{2}+1}\left({-\theta^1, -\theta^0};{t_1 \atop t_2};{q_1,p_1 \atop q_2,p_2}\right)\\
&=
t_2H_{\mathrm{III}(D_7)}\left(1-\theta^0; t_2; q_2,p_2\right)
-t_2p_1p_2-\frac{t_2}{t_1}p_1q_1+\frac{q_1q_2}{t_1}(p_1q_1+\theta^1). \nonumber
\end{align}

\subsection{Singularity pattern $\frac52+1$}
\subsubsection{Spectral type $(((1)(1)))_2, 31$}
The Riemann scheme is given by
\[
\left(
\begin{array}{cc}
  x=0  & x=\infty \, \left( \frac12 \right) \\
\begin{array}{c} 0 \\ 0 \\ 0 \\ \theta^{0} \end{array} &
\overbrace{\begin{array}{cccc}
1   & 0 & -t_1/2 & \theta^\infty_1/2\\
1   & 0 & -t_2/2 & \theta^\infty_2/2\\
-1 & 0 &   t_1/2 & \theta^\infty_1/2\\
-1 & 0 &   t_2/2 & \theta^\infty_2/2
      \end{array}}
\end{array}
\right),
\]
and the Fuchs-Hukuhara relation is written as
$\theta^0+\theta^\infty_1+\theta^\infty_2=0$.

The Lax pair is expressed as
\begin{equation}
\left\{
\begin{aligned}
\frac{\partial Y}{\partial x}&=
\left(
A_0x+A_1+\frac{A_2}{x}
\right)Y, \\
\frac{\partial Y}{\partial t_i}&=(-N_ix+B_i)Y \quad (i=1,2).
\end{aligned}
\right.
\end{equation}
Here
{\allowdisplaybreaks
\begin{align*}
A_0&=
\begin{pmatrix}
0 & 0 & 1 & 0 \\
0 & 0 & 0 & 1 \\
0 & 0 & 0 & 0 \\
0 & 0 & 0 & 0
\end{pmatrix}, \quad
N_1=
\begin{pmatrix}
0 & 0 & 1 & 0 \\
0 & 0 & 0 & 0 \\
0 & 0 & 0 & 0 \\
0 & 0 & 0 & 0
\end{pmatrix}, \quad
N_2=
\begin{pmatrix}
0 & 0 & 0 & 0 \\
0 & 0 & 0 & 1 \\
0 & 0 & 0 & 0 \\
0 & 0 & 0 & 0
\end{pmatrix}, \\
A_1&=
\begin{pmatrix}
0 & 0 & p_1-t_1 & p_1\\
0 & 0 & p_2 & p_2-t_2\\
1 & 0 & 0 & 0 \\
0 & 1 & 0 & 0
\end{pmatrix}, \quad
A_2=
\begin{pmatrix}
p_1q_1-\theta^\infty_1 \\
p_2q_2-\theta^\infty_2 \\
-p_1 \\
-p_2
\end{pmatrix}
\begin{pmatrix}
1 & 1 & q_1 & q_2
\end{pmatrix},\\
B_1&=
\begin{pmatrix}
\frac{p_1(q_1-q_2)-\theta^\infty_1}{t_1-t_2} & \frac{p_1(q_1-q_2)-\theta^\infty_1}{t_1-t_2} & -2p_1+t_1 & -p_1 \\
\frac{p_2(q_2-q_1)-\theta^\infty_2}{t_1-t_2} & \frac{p_2(q_2-q_1)-\theta^\infty_2}{t_1-t_2}-q_1 & -p_2 & 0 \\
-1 & 0 & \frac{p_1(q_1-q_2)-\theta^\infty_1}{t_1-t_2} & \frac{p_1(q_1-q_2)-\theta^\infty_1}{t_1-t_2} \\
0 & 0 & \frac{p_2(q_2-q_1)-\theta^\infty_2}{t_1-t_2} & \frac{p_2(q_2-q_1)-\theta^\infty_2}{t_1-t_2}-q_1
\end{pmatrix},\\
B_2&=
\begin{pmatrix}
\frac{p_1(q_1-q_2)-\theta^\infty_1}{t_2-t_1}-q_2 & \frac{p_1(q_1-q_2)-\theta^\infty_1}{t_2-t_1} & 0 & -p_1 \\
\frac{p_2(q_2-q_1)-\theta^\infty_2}{t_2-t_1} & \frac{p_2(q_2-q_1)-\theta^\infty_2}{t_2-t_1} & -p_2 & -2p_2+t_2 \\
0 & 0 & \frac{p_1(q_1-q_2)-\theta^\infty_1}{t_2-t_1}-q_2 & \frac{p_1(q_1-q_2)-\theta^\infty_1}{t_2-t_1} \\
0 & -1 & \frac{p_2(q_2-q_1)-\theta^\infty_2}{t_2-t_1} & \frac{p_2(q_2-q_1)-\theta^\infty_2}{t_2-t_1}
\end{pmatrix}.
\end{align*}
}

The Hamiltonians are given by
\begin{align}
&H_{\mathrm{Gar},t_1}^{\frac52+1+1}\left({\theta^\infty_1, \theta^\infty_2};{t_1 \atop t_2};{q_1,p_1\atop q_2,p_2}\right)\\
&\quad=H_\mathrm{II}\left(-\theta^\infty_1; t_1; q_1,p_1\right)
+p_1p_2+\frac{1}{t_1-t_2}(p_1(q_1-q_2)-\theta^\infty_1)(p_2(q_2-q_1)-\theta^\infty_2),\nonumber\\
&H_{\mathrm{Gar},t_2}^{\frac52+1+1}\left({\theta^\infty_1, \theta^\infty_2};{t_1 \atop t_2};{q_1,p_1\atop q_2,p_2}\right)\\
&\quad=H_\mathrm{II}\left(-\theta^\infty_2; t_2; q_2,p_2\right)
+p_1p_2+\frac{1}{t_2-t_1}(p_1(q_1-q_2)-\theta^\infty_1)(p_2(q_2-q_1)-\theta^\infty_2)\nonumber.
\end{align}

\subsubsection{Spectral type $((((((1))))))_4, 31$}
The Riemann scheme is given by
\[
\left(
\begin{array}{cc}
  x=0  & x=\infty \, \left( \frac14 \right) \\
\begin{array}{c} 0 \\ 0 \\ 0 \\ \theta^{0} \end{array} &
\overbrace{\begin{array}{ccccccc}
1 & 0 & 0 & 0 & -t_1/2 & \sqrt{t_2}/2 & \theta^\infty_1/4\\
-1 & 0 & 0 & 0 & t_1/2 & \sqrt{-1}\sqrt{t_2}/2 & \theta^\infty_1/4\\
1 & 0 & 0 & 0 & -t_1/2 & -\sqrt{t_2}/2 & \theta^\infty_1/4\\
-1 & 0 & 0 & 0 & t_1/2 & -\sqrt{-1}\sqrt{t_2}/2 & \theta^\infty_1/4
      \end{array}}
\end{array}
\right),
\]
and the Fuchs-Hukuhara relation is written as
$\theta^0+\theta^\infty_1=0$.

The Lax pair is expressed as
\begin{equation}
\left\{
\begin{aligned}
\frac{\partial Y}{\partial x}&=
\left(
A_0x+A_1+\frac{A_2}{x}
\right)Y, \\
\frac{\partial Y}{\partial t_i}&=(-N_ix+B_i)Y \quad (i=1,2).
\end{aligned}
\right.
\end{equation}
Here
{\allowdisplaybreaks
\begin{align*}
A_0&=
\begin{pmatrix}
0 & 0 & 1 & 0 \\
0 & 0 & 0 & 1 \\
0 & 0 & 0 & 0 \\
0 & 0 & 0 & 0
\end{pmatrix}, \quad
N_1=
\begin{pmatrix}
0 & 0 & 1 & 0 \\
0 & 0 & 0 & 1 \\
0 & 0 & 0 & 0 \\
0 & 0 & 0 & 0
\end{pmatrix}, \quad
N_2=
\begin{pmatrix}
0 & 0 & 0 & 1 \\
0 & 0 & 0 & 0 \\
0 & 0 & 0 & 0 \\
0 & 0 & 0 & 0
\end{pmatrix}, \\
A_1&=
\begin{pmatrix}
q_1 & q_2 & p_1-{q_1}^2-t_1 & -2q_1q_2-t_2 \\
0 & q_1 & p_2 & -{q_1}^2-t_1 \\
1 & 0 & -q_1 & -q_2 \\
0 & 1 & 0 & -q_1
\end{pmatrix}, \quad
A_2=
\begin{pmatrix}
\theta^0 \\
1 \\
-p_1 \\
-p_2
\end{pmatrix}
\begin{pmatrix}
1 & 0 & 0 & 0
\end{pmatrix},\\
B_1&=
\begin{pmatrix}
0 & 0 & 0 & 0 \\
0 & 0 & -2p_2 & 2p_1 \\
-1 & 0 & 2q_1 & 2q_2 \\
0 & -1 & 0 & 2q_1
\end{pmatrix}, \quad
B_2=\frac{1}{t_2}
\begin{pmatrix}
0 & 0 & 0 & 0 \\
1 & -2p_2q_2-\theta^0 &0 & -2q_2 \\
0 & -t_2 & 0 & 2t_2q_1 \\
0 & 0 & 1 & -2p_2q_2-\theta^0
\end{pmatrix}.
\end{align*}
}

The Hamiltonians are given by
\begin{align}
H^{\frac52+2}_{\mathrm{Gar},t_1}\left({\theta^0};{t_1 \atop t_2};{q_1,p_1 \atop q_2,p_2}\right)
&=
H_{\mathrm{II}}\left({\theta^0};t_1;q_1,p_1\right)
-2p_2q_2q_1-t_2p_2-q_2,\\
t_2H^{\frac52+2}_{\mathrm{Gar},t_2}\left({\theta^0};{t_1 \atop t_2};{q_1,p_1 \atop q_2,p_2}\right)
&={p_2}^2{q_2}^2+\theta^0 p_2q_2+t_2p_2(p_1-{q_1}^2-t_1)-p_1q_2-t_2q_1.
\end{align}

\subsection{Singularity pattern $2+\frac32$}
\subsubsection{Spectral type $(111)(1), (2)_2$}
The Riemann scheme is given by
\[
\left(
\begin{array}{cc}
  x=0 \, \left( \frac12 \right) & x=\infty \\
\overbrace{
\begin{array}{cc}
1 & 0 \\
1 & 0 \\
-1 & 0 \\
-1 & 0
\end{array}} &
\overbrace{
\begin{array}{cc}
-t & \theta^\infty_1\\
0 & \theta^\infty_2\\
0 & \theta^\infty_3\\
0 & \theta^\infty_4
\end{array}}
\end{array}
\right) ,
\]
and the Fuchs-Hukuhara relation is written as
$\theta_1^\infty+\theta_2^\infty+\theta_3^\infty+\theta_4^\infty =0$.

The Lax pair is expressed as
\begin{equation}
\left\{
\begin{aligned}
\frac{\partial Y}{\partial x}&=
\left(
\frac{A_2}{x^2}+\frac{A_1}{x}+A_0
\right)Y, \\
\frac{\partial Y}{\partial t}&=(E_1 x+B_0)Y.
\end{aligned}
\right.
\end{equation}
Here
{\allowdisplaybreaks
\begin{align*}
A_{\xi}&=
U^{-1}\hat{A}_{\xi}U, \quad
A_0=tE_{1}, \quad
B_0=U^{-1}\hat{B}_0U, \quad
 U=\mathrm{diag}(1,u,v,w), \quad E_1=\mathrm{diag}(1,0,0,0), \\
\hat{A}_1&=
\begin{pmatrix}
-\theta^\infty_1 & a_{12} & (1-f_1)q_1-q_2 & -q_1 \\
-p_1q_1+\theta^\infty_2+\theta^\infty_4 & -\theta^\infty_2 & 0 & 0 \\
a_{31} & 0 & -\theta^\infty_3 & 0 \\
a_{41} & 0 & 0 & -\theta^\infty_4
\end{pmatrix},\\
\hat{A}_2&=
\begin{pmatrix}
I_2 \\ B_2
\end{pmatrix}
\begin{pmatrix}
C_2B_2 & -C_2
\end{pmatrix}, \quad 
B_2=\begin{pmatrix}
p_2 & f_2 \\
p_1+(1-f_1)p_2 & f_3
\end{pmatrix}, \\
C_2&=\begin{pmatrix}
f_1 & 1 \\
f_1-1 & 1
\end{pmatrix}, \quad
\hat{B}_0=\frac{1}{t}
\begin{pmatrix}
0 & (\hat{A}_1)_{12} & (\hat{A}_1)_{13} & (\hat{A}_1)_{14}\\
(\hat{A}_1)_{21} & 0 & 0 & 0\\
(\hat{A}_1)_{31} & 0 & 0 & 0 \\
(\hat{A}_1)_{41} & 0 & 0 & 0
\end{pmatrix},
\end{align*}
}
where
\begin{align*}
&(\theta^\infty_3-\theta^\infty_4)f_1=-p_1(q_1-q_2)+\theta^\infty_2+\theta^\infty_3,\\
&(\theta^\infty_2-\theta^\infty_3)f_2=p_2^2(q_2-q_1)-(\theta^\infty_2+\theta^\infty_3)p_2+1,\\
&(\theta^\infty_4-\theta^\infty_2)f_3=((1-f_1)p_2+p_1)((q_1-q_2)p_2+\theta^\infty_2+\theta^\infty_3)+f_1,\\
&a_{12}=f_3q_1+p_2(q_2-q_1)+f_2((f_1-1)q_1+q_2)-\theta^\infty_2-\theta^\infty_3,\\
&a_{31}=(\theta^\infty_3-\theta^\infty_1)p_2-p_2^2q_2-p_1q_1(f_2+p_2)+(\theta^\infty_2+\theta^\infty_4)f_2-1,\\
&a_{41}=(p_1q_1+p_2q_2+\theta^\infty_1-\theta^\infty_4)((f_1-1)p_2-p_1)+ f_3(-p_1q_1+\theta^\infty_2+\theta^\infty_4)+f_1-1.
\end{align*}

The Hamiltonian is given by
\begin{align}
&tH_\mathrm{KSs}^{2+\frac32}\left( {-2\theta^\infty_4, -2\theta^\infty_3 \atop \theta^\infty_1+\theta^\infty_3}; t; {q_1,p_1\atop q_2,p_2} \right)\\
&=tH_{\mathrm{III}(D_7)}\left( -2\theta^\infty_4; t ; q_1, p_1 \right)
+tH_{\mathrm{III}(D_7)}\left( -2\theta^\infty_3; t; q_2, p_2 \right)
+2p_2q_1(p_1q_1+\theta^\infty_1+\theta^\infty_3)-q_1.\nonumber
\end{align}

The gauge parameters $u$, $v$, $w$ satisfy
\begin{align}\label{eq:gauge_(1)_2,(2)(2)}
\frac{t}{u}\frac{du}{dt}&=-2 p_2q_1+\theta^\infty_1-\theta^\infty_2,\quad
\frac{t}{v}\frac{dv}{dt}=-2p_2q_2+\theta^\infty_1+\theta^\infty_3, \\
\frac{t}{w}\frac{dw}{dt}&=-2(p_1+p_2)q_1-\theta^\infty_2-\theta^\infty_3.\nonumber
\end{align}

\subsubsection{Spectral type $(1)_2 11, (2)(2)$}
The Riemann scheme is given by
\[
\left(
\begin{array}{cc}
  x=0  & x=\infty \, \left( \frac12 \right) \\
\overbrace{
\begin{array}{cc}
0 & 0 \\
0 & 0 \\
1 & \theta^0 \\
1 & \theta^0
\end{array}} &
\overbrace{\begin{array}{cc}
0         & \theta^\infty_2 \\
0         & \theta^\infty_3 \\
\sqrt{t} & \theta^\infty_1/2\\
-\sqrt{t}& \theta^\infty_1/2
      \end{array}}
\end{array}
\right) ,
\]
and the Fuchs-Hukuhara relation is written as
$2\theta^0+\theta^\infty_1+\theta^\infty_2+\theta^\infty_3=0$.

The Lax pair is expressed as
\begin{equation}
\left\{
\begin{aligned}
\frac{\partial Y}{\partial x}&=
\left(
\frac{A_2}{x^2}+\frac{A_1}{x}+A_0
\right)Y, \\
\frac{\partial Y}{\partial t}&=\left( \frac1t A_0 x+B_0 \right)Y.
\end{aligned}
\right.
\end{equation}
Here
{\allowdisplaybreaks
\begin{align*}
A_0&=
\begin{pmatrix}
0 & t & 0 & 0 \\
0 & 0 & 0 & 0 \\
0 & 0 & 0 & 0\\
0 & 0 & 0 & 0
\end{pmatrix}, \quad
A_2=
\begin{pmatrix}
0 & 0 & 0 & 0 \\
p_1+p_2 & 1 & 0 & p_1 \\
p_1 & 0 & 1 & p_1 \\
0 & 0 & 0 & 0
\end{pmatrix},\\
A_1&=
{\footnotesize
\begin{pmatrix}
p_1q_1+p_2q_2-\theta^0-\theta^\infty_1 & q_2 & q_1-q_2-\frac{\theta^0+\theta^\infty_1}{p_1} & p_1q_1-\theta^0-\theta^\infty_1 \\
1 &  -p_1q_1-p_2q_2+\theta^0 & p_2(q_2-q_1)+(\theta^0+\theta^\infty_1)\frac{p_2}{p_1}-\theta^0-\theta^\infty_2 & 0 \\
0 & -p_1q_1 & -\theta^\infty_2 & 0 \\
p_2(q_1-q_2)-\theta^\infty_3 & q_1-q_2 & q_2-q_1-\frac{\theta^\infty_3}{p_1} & -\theta^\infty_3
\end{pmatrix}
},\\
B_0&=\frac{1}{t}
{\footnotesize
\begin{pmatrix}
p_1q_1+p_2q_2-\theta^0 & 0 & 0 & p_1q_1-\theta^0-\theta^\infty_1 \\
1 & -p_1q_1-p_2q_2+\theta^0+\theta^\infty_1 & p_2(q_2-q_1)+(\theta^0+\theta^\infty_1)\frac{p_2}{p_1}-\theta^0-\theta^\infty_2 & 0\\
0 & -p_1q_1 & -2p_2q_1+2(\theta^0+\theta^\infty_1)\frac{p_2}{p_1}+\theta^\infty_1-\theta^\infty_2 & 0 \\
p_2(q_1-q_2)-\theta^\infty_3 & 0 & 0 & 2p_2q_1+\theta^\infty_1-\theta^\infty_3
\end{pmatrix}
}.
\end{align*}
}

The Hamiltonian is given by
\begin{align}
&tH_\mathrm{KSs}^{2+\frac32}\left( {\theta^0+2\theta^\infty_2, -\theta^0 \atop -\theta^0-\theta^\infty_1}; t; {q_1, p_1\atop q_2, p_2} \right)\\
&=tH_{\mathrm{III}(D_7)}\left( \theta^0+2\theta^\infty_2; t ;q_1,p_1 \right)
+tH_{\mathrm{III}(D_7)}\left( -\theta^0; t; q_2, p_2 \right)
+2p_2q_1(p_1q_1-\theta^0-\theta^\infty_1)-q_1.\nonumber
\end{align}

\subsection{Singularity pattern $\frac43+1+1$}
\subsubsection{Spectral type $(1)_3 1, 22, 22$}
The Riemann scheme is given by
\[
\left(
\begin{array}{ccc}
  x=0 & x=1 & x=\infty \, \left( \frac13 \right) \\
\begin{array}{c} 0 \\ 0 \\ \theta^0 \\ \theta^0 \end{array} &
\begin{array}{c} 0 \\ 0 \\ \theta^1 \\ \theta^1 \end{array} &
\overbrace{
\begin{array}{cc}
0 & \theta^\infty_2 \\
-t^{\frac13}               & \theta^\infty_1/3 \\
-\omega t^{\frac13}    & \theta^\infty_1/3 \\
-\omega^2 t^{\frac13} & \theta^\infty_1/3
\end{array}}
\end{array}
\right) ,
\]
and the Fuchs-Hukuhara relation is written as
$2\theta^0+2\theta^1+\theta^\infty_1+\theta^\infty_2=0$.
Here $\omega=e^{\frac{2}{3}\pi i}$ is a cube root of unity. 

The Lax pair is expressed as
\begin{equation}
\left\{
\begin{aligned}
\frac{\partial Y}{\partial x}&=
\left(
\frac{A_0}{x}+\frac{A_1}{x-1}+A_\infty
\right)Y, \\
\frac{\partial Y}{\partial t}&=\left( \frac1t A_\infty x+B \right)Y.
\end{aligned}
\right.
\end{equation}
Here
{\allowdisplaybreaks
\begin{align*}
A_\infty&=
\begin{pmatrix}
0 & 1 & 0 & 0 \\
0 & 0 & 0 & 0 \\
0 & 0 & 0 & 0 \\
0 & 0 & 0 & 0
\end{pmatrix}, \quad
A_0=
\begin{pmatrix}
0 & 0 & 0 & 0 \\
-q_2 & \theta^0 & 0 & q_1(p_1q_1-\theta^1)-\theta^0 q_2+t \\
1 & 0 & \theta^0 & 0 \\
0 & 0 & 0 & 0
\end{pmatrix}, \\
A_1&=
\begin{pmatrix}
B_1 \\
I_2
\end{pmatrix}
\begin{pmatrix}
C_1 & \theta^1 I_2-C_1 B_1
\end{pmatrix}, \\
B_1&=
\begin{pmatrix}
p_2(q_2-q_1)-\theta^0-\theta^1-\theta^\infty_2 & -p_2q_2+\theta^1+\theta^\infty_2 \\
q_1-q_2 & q_2
\end{pmatrix}, \quad
C_1=
\begin{pmatrix}
0 & p_1 \\
1 & p_1+p_2
\end{pmatrix},\\
B&=\frac{1}{t}
{\small
\begin{pmatrix}
-p_2q_2+\theta^1+\theta^\infty_2 & 0 & 0 & B_{14} \\
0 & p_1q_1+p_2q_2+\theta^0+1 & B_{23} & B_{24} \\
1 & -p_2 & -2p_2q_2-\theta^0-\theta^1-\theta^\infty_1 & -p_1q_1+2p_2q_2-\theta^\infty_2 \\
1 & 0 & 0 & -p_1q_1+\theta^0+\theta^1 \\
\end{pmatrix}},
\end{align*}
where
\begin{align*}
B_{14}&=-p_2((p_1q_1-\theta^1)q_1+\theta^0 q_2+t)+\theta^0(p_1q_1+\theta^\infty_2),\\
B_{23}&=p_1q_1(q_2-q_1)+\theta^1 q_1+(\theta^0+\theta^\infty_2)q_2,\\
B_{24}&=p_1q_1(q_1-q_2)-\theta^1 q_1-(\theta^0+\theta^\infty_2)q_2+t.
\end{align*}
}

The Hamiltonian is given by
\begin{align}
&tH_\mathrm{KSs}^{2+\frac32}\left( {\theta^0-\theta^1+1, -\theta^0-\theta^1-2\theta^\infty_2+1 \atop -\theta^1}; t; {q_1, p_1 \atop q_2, p_2} \right)\\
&=tH_{\mathrm{III}(D_7)}\left( \theta^0-\theta^1+1; t; q_1, p_1 \right)
+tH_{\mathrm{III}(D_7)}\left( -\theta^0-\theta^1-2\theta^\infty_2+1; t; q_2, p_2 \right)\nonumber\\
&\quad+2p_2q_1(p_1q_1-\theta^1)-q_1.\nonumber
\end{align}

\subsection{Singularity pattern $\frac32+\frac32$}
\subsubsection{Spectral type $(1)_2 11, (2)_2$}
The Riemann scheme is given by
\[
\left(
\begin{array}{cc}
  x=0 \, \left( \frac12 \right) & x=\infty \, \left( \frac12 \right) \\
\overbrace{
\begin{array}{cc}
1 & 0 \\
1 & 0 \\
-1 & 0 \\
-1 & 0
\end{array}} &
\overbrace{\begin{array}{cc}
0         & \theta^\infty_2 \\
0         & \theta^\infty_3 \\
\sqrt{t} & \theta^\infty_1/2\\
-\sqrt{t}& \theta^\infty_1/2
      \end{array}}
\end{array}
\right) ,
\]
and the Fuchs-Hukuhara relation is written as
$\theta^\infty_1+\theta^\infty_2+\theta^\infty_3=0$.

The Lax pair is expressed as
\begin{equation}
\left\{
\begin{aligned}
\frac{\partial Y}{\partial x}&=
\left(
\frac{A_2}{x^2}+\frac{A_1}{x}+A_0
\right)Y, \\
\frac{\partial Y}{\partial t}&=\left( \frac1t A_0 x+B_0 \right)Y.
\end{aligned}
\right.
\end{equation}
Here
{\allowdisplaybreaks
\begin{align*}
A_0&=
\begin{pmatrix}
0 & 1 & 0 & 0 \\
0 & 0 & 0 & 0 \\
0 & 0 & 0 & 0\\
0 & 0 & 0 & 0
\end{pmatrix}, \quad
A_2=
\begin{pmatrix}
0   & 0 & 0 & 0 \\
q_2 & 0 & 0 & q_1q_2-t \\
1   & 0 & 0 & q_1 \\
0   & 0 & 0 & 0
\end{pmatrix},\\
A_1&=
{\small
\begin{pmatrix}
p_1q_1-p_2q_2+\theta^\infty_2 & q_1/t & 1-\frac{q_1q_2}{t} & -p_2(q_1q_2-t)-\theta^\infty_1q_1 \\
t &  -p_1q_1+p_2q_2+\theta^\infty_3 & p_1(q_1q_2-t)+\theta^\infty_1q_2 & 0 \\
0 & p_2 & -\theta^\infty_2 & 0 \\
-p_1 & -1/t & q_2/t & -\theta^\infty_3
\end{pmatrix}
},\\
B_0&=\frac{1}{t}
{\small
\begin{pmatrix}
p_1q_1-p_2q_2-\theta^\infty_3 & 0 & 0 & -p_2(q_1q_2-t)-\theta^\infty_1q_1 \\
t & -p_1q_1+p_2q_2-\theta^\infty_2+1 & p_1(q_1q_2-t)+\theta^\infty_1q_2 & 0\\
0 & p_2 & -2p_2q_2-\theta^\infty_3 & 0 \\
-p_1 & 0 & 0 & -2p_1q_1-\theta^\infty_2
\end{pmatrix}
}.
\end{align*}
}

The Hamiltonian is given by
\begin{align}
&tH_\mathrm{KSs}^{2+\frac43}\left( {-2\theta^\infty_3, -2\theta^\infty_2+1}; t; {q_1, p_1 \atop q_2, p_2} \right)\\
&=tH_{\mathrm{III}(D_7)}\left( -2\theta^\infty_3; t ;q_1, p_1 \right)
+tH_{\mathrm{III}(D_7)}\left( -2\theta^\infty_2+1; t; q_2, p_2 \right)
-t(2p_1p_2+p_1+p_2).\nonumber
\end{align}

\subsection{Singularity pattern $2+\frac43$}
\subsubsection{Spectral type $(1)_3 1,(2)(2)$}
The Riemann scheme is given by
\[
\left(
\begin{array}{cc}
  x=0  & x=\infty \, \left( \frac13 \right) \\
\overbrace{
\begin{array}{cc}
0 & 0 \\
0 & 0 \\
1 & \theta^0 \\
1 & \theta^0
\end{array}} &
\overbrace{\begin{array}{cc}
0         & \theta^\infty_2 \\
t^{\frac13}               & \theta^\infty_1/3 \\
\omega t^{\frac13}    & \theta^\infty_1/3 \\
\omega^2 t^{\frac13} & \theta^\infty_1/3
      \end{array}}
\end{array}
\right) ,
\]
and the Fuchs-Hukuhara relation is written as
$2\theta^0+\theta^\infty_1+\theta^\infty_2=0$.

The Lax pair is expressed as
\begin{equation}
\left\{
\begin{aligned}
\frac{\partial Y}{\partial x}&=
\left(
\frac{A_2}{x^2}+\frac{A_1}{x}+A_0
\right)Y, \\
\frac{\partial Y}{\partial t}&=\left( \frac1t A_0 x+B_0 \right)Y.
\end{aligned}
\right.
\end{equation}
Here
{\allowdisplaybreaks
\begin{align*}
A_0&=
\begin{pmatrix}
0 & 1 & 0 & 0 \\
0 & 0 & 0 & 0 \\
0 & 0 & 0 & 0\\
0 & 0 & 0 & 0
\end{pmatrix}, \quad
A_2=
\begin{pmatrix}
0 & 0 & 0 & 0 \\
q_1 & 1 & 0 & p_2 \\
1 & 0 & 1 & 0 \\
0 & 0 & 0 & 0
\end{pmatrix},\\
A_1&=
{\small
\begin{pmatrix}
-p_1q_1+\theta^0+\theta^\infty_2 & -p_1 & \theta^0+\theta^\infty_2 & -p_1p_2 \\
0 &  p_1q_1-p_2q_2+\theta^0 & p_2(q_1q_2-t)-(\theta^0+\theta^\infty_2)q_1 & 1 \\
0 & p_1 & -\theta^\infty_2 & 0 \\
t & q_2 & -q_1q_2+t & p_2q_2
\end{pmatrix}
},\\
B_0&=\frac{1}{t}
{\small
\begin{pmatrix}
-p_1q_1+\theta^0+\theta^\infty_2 & 0 & 0 & p_1p_2 \\
0 & p_1q_1-p_2q_2+\theta^0+1 & p_2(q_1q_2-t)-(\theta^0+\theta^\infty_2)q_1 & 1\\
0 & p_1 & -2p_1q_1+\theta^0+\theta^\infty_2 & 0 \\
t & 0 & 0 & p_2q_2+1
\end{pmatrix}
}.
\end{align*}
}
The Hamiltonian is given by
\begin{align}
&tH_\mathrm{KSs}^{2+\frac43}\left( -\theta^0-2\theta^\infty_2+1, -\theta^0 ;t; {q_1,p_1 \atop q_2,p_2} \right)\\
&=tH_{\mathrm{III}(D_7)}\left( -\theta^0-2\theta^\infty_2+1; t; q_1, p_1 \right)
+tH_{\mathrm{III}(D_7)}\left( -\theta^0; t; q_2, p_2 \right)
-t(2p_1p_2+p_1+p_2).\nonumber
\end{align}

\subsection{Singularity pattern $\frac54+1+1$}
\subsubsection{Spectral type $(1)_4, 22, 22$}
The Riemann scheme is given by
\[
\left(
\begin{array}{ccc}
  x=0 & x=1 & x=\infty \, \left( \frac14 \right) \\
\begin{array}{c} 0 \\ 0 \\ \theta^0 \\ \theta^0 \end{array} &
\begin{array}{c} 0 \\ 0 \\ \theta^1 \\ \theta^1 \end{array} &
\overbrace{
\begin{array}{cc}
t^{\frac14} & \theta^\infty_1/4 \\
\sqrt{-1}t^{\frac14}     & \theta^\infty_1/4 \\
-t^{\frac14}    & \theta^\infty_1/4 \\
-\sqrt{-1}t^{\frac14} & \theta^\infty_1/4
\end{array}}
\end{array}
\right) ,
\]
and the Fuchs-Hukuhara relation is written as
$2\theta^0+2\theta^1+\theta^\infty_1=0$.

The Lax pair is expressed as
\begin{equation}
\left\{
\begin{aligned}
\frac{\partial Y}{\partial x}&=
\left(
\frac{A_0}{x}+\frac{A_1}{x-1}+A_\infty
\right)Y, \\
\frac{\partial Y}{\partial t}&=\left( \frac1t A_\infty x+B \right)Y.
\end{aligned}
\right.
\end{equation}
Here
{\allowdisplaybreaks
\begin{align*}
A_\infty&=
\begin{pmatrix}
0 & 1 & 0 & 0 \\
0 & 0 & 0 & 0 \\
0 & 0 & 0 & 0 \\
0 & 0 & 0 & 0
\end{pmatrix}, \quad
A_0=
\begin{pmatrix}
-p_1q_1 & p_1(p_1q_1+\theta^0) & -p_1p_2 & 0 \\
-q_1 & p_1q_1+\theta^0 & -p_2 & 0 \\
0 & 0 & 0 & 0 \\
t p_2 & -t p_1p_2 & p_2(p_2q_2+\theta^1-\theta^0)+1 & \theta^0
\end{pmatrix}, \\
A_1&=
\begin{pmatrix}
0 & 0 & 0 & 0 \\
q_1 & \theta^1 & 0 & 1 \\
t & 0 & \theta^1 & q_2 \\
0 & 0 & 0 & 0
\end{pmatrix},\\
B&=\frac{1}{t}
\begin{pmatrix}
-p_1q_1 & -1 & 0 & p_1 \\
0 & p_1q_1+\theta^0+\theta^1+1 & -p_2 & 1 \\
t & -t p_1 & p_2q_2+\theta^1+1 & 0 \\
t p_2 & 0 & 0 & -p_2q_2+\theta^0+1 \\
\end{pmatrix}.
\end{align*}
}

The Hamiltonian is given by
\begin{align}
&tH_\mathrm{KSs}^{2+\frac43}\left( {\theta^0+\theta^1+1, \theta^1-\theta^0}; t; {q_1, p_1 \atop q_2, p_2} \right)\\
&=tH_{\mathrm{III}(D_7)}\left( \theta^0+\theta^1+1; t; q_1, p_1 \right)
+tH_{\mathrm{III}(D_7)}\left( \theta^1-\theta^0; t; q_2, p_2 \right)-t(2p_1p_2+p_1+p_2).\nonumber
\end{align}

\subsection{Singularity pattern $\frac32+\frac43$}
\subsubsection{Spectral type $(1)_3 1, (2)_2$}
The Riemann scheme is given by
\[
\left(
\begin{array}{cc}
  x=0 \, \left( \frac12 \right) & x=\infty \, \left( \frac13 \right) \\
\overbrace{
\begin{array}{cc}
1 & 0 \\
1 & 0 \\
-1 & 0 \\
-1 & 0
\end{array}} &
\overbrace{\begin{array}{cc}
0                            & \theta^\infty_2 \\
-t^{\frac13}               & \theta^\infty_1/3 \\
-\omega t^{\frac13}    & \theta^\infty_1/3 \\
-\omega^2 t^{\frac13} & \theta^\infty_1/3
      \end{array}}
\end{array}
\right) ,
\]
and the Fuchs-Hukuhara relation is written as
$\theta^\infty_1+\theta^\infty_2=0$.

The Lax pair is expressed as
\begin{equation}
\left\{
\begin{aligned}
\frac{\partial Y}{\partial x}&=
\left(
\frac{A_2}{x^2}+\frac{A_1}{x}+A_0
\right)Y, \\
\frac{\partial Y}{\partial t}&=\left( \frac1t A_0 x+B_0 \right)Y.
\end{aligned}
\right.
\end{equation}
Here
{\allowdisplaybreaks
\begin{align*}
A_0&=
\begin{pmatrix}
0 & 1 & 0 & 0 \\
0 & 0 & 0 & 0 \\
0 & 0 & 0 & 0\\
0 & 0 & 0 & 0
\end{pmatrix}, \quad
A_2=
\begin{pmatrix}
0   & 0 & 0 & 0 \\
q_2 & 0 & 0 & 1 \\
1   & 0 & 0 & 0 \\
0   & 0 & 0 & 0
\end{pmatrix},\\
A_1&=
\begin{pmatrix}
-p_2q_2+\theta^\infty_2+\frac12 & 0 & 1 & -p_2 \\
0 &  -p_1q_1+p_2q_2+\theta^\infty_2 & -\frac{t}{q_1}+q_2\left( p_1q_1-2\theta^\infty_2-\frac12 \right) & q_1 \\
0 & p_2 & -\theta^\infty_2-\frac12 & 0 \\
t/q_1 & 1 & -q_2 & p_1q_1-\theta^\infty_2
\end{pmatrix},\\
B_0&=\frac{1}{t}
\begin{pmatrix}
-p_2q_2+\frac12 & 0 & 0 & p_2 \\
0 & -p_1q_1+p_2q_2+1 & -\frac{t}{q_1}+q_2\left( p_1q_1-2\theta^\infty_2-\frac12\right) & q_1\\
0 & p_2 & -2p_2q_2+\frac12 & 0 \\
t/q_1 & 0 & 0 & -p_1q_1+1
\end{pmatrix}.
\end{align*}
}

The Hamiltonian is given by
\begin{align}
&tH_\mathrm{KSs}^{2+\frac54}\left( -2\theta^\infty_2; t; {q_1, p_1 \atop q_2, p_2} \right)  \\
&\quad=
tH_{\mathrm{III}(D_7)}\left( -2\theta^\infty_2; t ;q_1, p_1 \right)
+tH_{\mathrm{III}(D_7)}\left( -2\theta^\infty_2; t; q_2, p_2 \right)
-t\left( 2\frac{p_2}{q_1}+p_1+p_2 \right).\nonumber
\end{align}

\subsection{Singularity pattern $2+\frac54$}
\subsubsection{Spectral type $(1)_4, (2)(2)$}
The Riemann scheme is given by
\[
\left(
\begin{array}{cc}
  x=0  & x=\infty \, \left( \frac14 \right) \\
\overbrace{
\begin{array}{cc}
0 & 0 \\
0 & 0 \\
1 & \theta^0 \\
1 & \theta^0
\end{array}} &
\overbrace{\begin{array}{cc}
t^{\frac14}              & \theta^\infty_1/4 \\
\sqrt{-1}t^{\frac14}   & \theta^\infty_1/4 \\
-t^{\frac14}            & \theta^\infty_1/4 \\
-\sqrt{-1}t^{\frac14} & \theta^\infty_1/4
      \end{array}}
\end{array}
\right) ,
\]
and the Fuchs-Hukuhara relation is written as
$2\theta^0+\theta^\infty_1=0$.

The Lax pair is expressed as
\begin{equation}
\left\{
\begin{aligned}
\frac{\partial Y}{\partial x}&=
\left(
\frac{A_2}{x^2}+\frac{A_1}{x}+A_0
\right)Y, \\
\frac{\partial Y}{\partial t}&=\left( \frac1t A_0 x+B_0 \right)Y.
\end{aligned}
\right.
\end{equation}
Here
{\allowdisplaybreaks
\begin{align*}
A_0&=
\begin{pmatrix}
0 & t & 0 & 0 \\
0 & 0 & 0 & 0 \\
0 & 0 & 0 & 0\\
0 & 0 & 0 & 0
\end{pmatrix}, \quad
A_1=
{\small
\begin{pmatrix}
-p_1q_1+\theta^0+\frac12 & \frac{t}{q_1}\left( \theta^0+\frac12 \right)-tp_1 & -t/q_1 & 0 \\
0 &  p_1q_1-\frac12 & 1 & -p_2 \\
p_2 & 0 & \theta^0 & p_2(p_2q_2-\theta^0)+1 \\
1 & t/q_1 & q_2 & 0
\end{pmatrix}
},\\
A_2&=
\begin{pmatrix}
0 & 0 & 0 & 0 \\
q_1/t & 1 & 0 & 0 \\
0 & 0 & 1 & 0 \\
0 & 0 & 0 & 0
\end{pmatrix}, \quad
B_0=\frac{1}{t}
{\small
\begin{pmatrix}
-p_1q_1+\theta^0+\frac12 & 0 & 0 & t p_2/q_1 \\
0 & p_1q_1-\frac12 & 1 & -p_2\\
p_2 & t p_2/q_1 & -p_2q_2+\theta^0 & 0 \\
1 & 0 & 0 & p_2q_2
\end{pmatrix}
}.
\end{align*}
}

The Hamiltonian is given by
\begin{align}
&tH_\mathrm{KSs}^{2+\frac54}\left( -\theta^0 ;t; {q_1,p_1\atop q_2,p_2} \right)=
tH_{\mathrm{III}(D_7)}\left( -\theta^0; t ;q_1, p_1 \right)
+tH_{\mathrm{III}(D_7)}\left( -\theta^0; t; q_2, p_2 \right)
-t\left( 2\frac{p_2}{q_1}+p_1+p_2 \right).
\end{align}

\subsection{Singularity pattern $\frac32+\frac54$}
\subsubsection{Spectral type $(1)_4, (2)_2$}
The Riemann scheme is given by
\[
\left(
\begin{array}{cc}
  x=0 \, \left( \frac12 \right) & x=\infty \, \left( \frac14 \right) \\
\overbrace{
\begin{array}{cc}
1 & 0 \\
1 & 0 \\
-1 & 0 \\
-1 & 0
\end{array}} &
\overbrace{\begin{array}{cc}
t^{\frac14}              & 0 \\
\sqrt{-1}t^{\frac14}   & 0 \\
-t^{\frac14}            & 0 \\
-\sqrt{-1}t^{\frac14} & 0
      \end{array}}
\end{array}
\right).
\]
The Lax pair is expressed as
\begin{equation}\label{eq:(1)_4,(2)_2}
\left\{
\begin{aligned}
\frac{\partial Y}{\partial x}&=
\left(
\frac{A_2}{x^2}+\frac{A_1}{x}+A_0
\right)Y, \\
\frac{\partial Y}{\partial t}&=\left( \frac1t A_0 x+B_0 \right)Y.
\end{aligned}
\right.
\end{equation}
Here
\begin{align*}
A_0&=
\begin{pmatrix}
0 & t & 0 & 0 \\
0 & 0 & 0 & 0 \\
0 & 0 & 0 & 0\\
0 & 0 & 0 & 0
\end{pmatrix}, \quad
A_1=
\begin{pmatrix}
p_2q_2+\frac12 & -q_2 & 0 & \frac{q_1q_2}{t} \\
0 & -p_2q_2-\frac12 & q_1/t & 0 \\
0 & -q_2 & p_1q_1 & -t/q_1 \\
-1 & 0 & 1 & -p_1q_1
\end{pmatrix},\\
A_2&=
\begin{pmatrix}
0 & 0 & 0 & 0 \\
-1/q_2 & 0 & 0 & 0 \\
0 & 0 & 0 & 1 \\
0 & 0 & 0 & 0
\end{pmatrix}, \quad
B_0=\frac{1}{t}
\begin{pmatrix}
p_2q_2+\frac12 &  0 & 0 & -\frac{q_1q_2}{t} \\
0 & -p_2q_2-\frac12 & q_1/t & 0 \\
0 & q_2 & -p_1q_1 & -t/q_1 \\
-1 & 0 & 0 & -p_1q_1
\end{pmatrix}.
\end{align*}

The Hamiltonian is given by
\begin{align}
&tH_\mathrm{KSs}^{\frac32+\frac54}\left( t; {q_1, p_1 \atop q_2, p_2} \right)=
tH_{\mathrm{III}(D_8)}\left(t ;q_1,p_1\right)+tH_{\mathrm{III}(D_8)}\left(t; q_2,p_2\right)
-2\frac{q_1q_2}{t}+q_1+q_2.
\end{align}

\section{Laplace transform}
In the degeneration scheme of the Sasano system,
there are some Hamiltonians which have more than one associated linear systems;
that is, $H^{D_5}_{\mathrm{Ss}}$, $H^{D_4}_{\mathrm{Ss}}$, $H_\mathrm{KSs}^{2+\frac32}$, $H_\mathrm{KSs}^{2+\frac43}$, and $H_\mathrm{KSs}^{2+\frac54}$.

Linear systems associated with the same Hamiltonian can be transformed into one another
by the Laplace transform (and a change of the independent variable if necessary).

The correspondences through the Laplace transform involving only unramified linear systems are
given in \cite{KNS}.
In this section we describe the correspondences involving ramified linear systems.

In the case of linear systems of the following form:
\begin{equation}
\frac{d}{dx}Y=\left[ B\left(xI_l-T\right)^{-1}C+S\right] Y
\end{equation}
where $B$ is an $m \times l$ matrix and $C$ is an $l \times m$ matrix,
its Laplace transform is the following~(\cite{B2, Hrd}).
\begin{equation}
\frac{d}{d\xi}\hat{Z}=-\left[ C\left(\xi I_m-S\right)^{-1}B+T\right]\hat{Z}.
\end{equation}
Through this transformation, we have the following correspondences of spectral types:
\begin{align*}
H^{D_4}_{\mathrm{Ss}}&:\  \stackrel{\infty}{(2)}_2, 31, 1111 \leftrightarrow \,
\stackrel{\infty}{(111)(1)}, (2)(2), \quad
\stackrel{\infty}{(1)_2 11}, 22, 22 \leftrightarrow \, \stackrel{\infty}{(2)(2)}, (111)(1), \\
H^{\frac32+1+1+1}_{\mathrm{Gar}}&:\  \stackrel{\infty}{(1)_2(1)_2}, 31, 22 \leftrightarrow \,
\stackrel{\infty}{(2)(1)}, (1)(1)(1), \\
H^{2+\frac32+1}_{\mathrm{Gar}}&:\  \stackrel{\infty}{((1))_4}, 31, 22 \leftrightarrow \,
\stackrel{\infty}{(2)(1)}, (1)_2(1), \\
H_\mathrm{KSs}^{2+\frac32}&:\ \stackrel{\infty}{(1)_3 1}, 22, 22 \leftrightarrow \,
\stackrel{\infty}{(2)(2)}, (1)_2 11, \quad
\stackrel{\infty}{(2)_2}, (111)(1) \leftrightarrow \,
\stackrel{\infty}{(1)_2 11}, (2)(2), \\
H_\mathrm{KSs}^{2+\frac43}&:\ \stackrel{\infty}{(1)_4}, 22, 22 \leftrightarrow \,
\stackrel{\infty}{(2)(2)}, (1)_3 1, \quad
\stackrel{\infty}{(2)_2}, (1)_2 11 \leftrightarrow \,
\stackrel{\infty}{(1)_3 1}, (2)(2), \\
H_\mathrm{KSs}^{2+\frac54}&:\ \stackrel{\infty}{(2)_2}, (1)_3 1 \leftrightarrow \,
\stackrel{\infty}{(1)_4}, (2)(2).
\end{align*}
Note that the symbol $\infty$ over each spectral type indicates which spectral type (of a singular point) corresponds to the singular point $x=\infty$.

Here the Hamiltonian $H^{\frac32+1+1+1}_{\mathrm{Gar}}$ is associated with only one linear system in
the degeneration scheme of the Sasano system.
However it has another linear system
which appears as the degeneration of 21,21,111,111.
They are related by this transformation.
The same holds for $H^{2+\frac32+1}_{\mathrm{Gar}}$.

Similarly, the linear system $\stackrel{\infty}{(((1)(1)))_2}, 31$ and $\stackrel{\infty}{((((((1))))))_4}, 31$,
which are associated with $H^{\frac52+1+1}_{\mathrm{Gar}}$ and $H^{\frac52+2}_{\mathrm{Gar}}$ respectively,
correspond to the system $\stackrel{\infty}{(((1)(1)))(((1)))}$ and $\stackrel{\infty}{(((1)))(1)_2}$
in the same manner as Section~4.3 in \cite{KNS} or Section~4 in \cite{K1}:
\begin{align*}
H^{\frac52+1+1}_{\mathrm{Gar}}&:\  \stackrel{\infty}{(((1)(1)))_2}, 31 \leftrightarrow \,
\stackrel{\infty}{(((1)(1)))(((1)))}, \\
H^{\frac52+2}_{\mathrm{Gar}}&:\  \stackrel{\infty}{((((((1))))))_4}, 31 \leftrightarrow \,
\stackrel{\infty}{(((1)))(1)_2}.
\end{align*}
Here the linear systems $\stackrel{\infty}{(((1)(1)))(((1)))}$ and $\stackrel{\infty}{(((1)))(1)_2}$ are
degenerated systems derived from the $21,21,111,111$-system, see~\cite{KNS, K3}.

\appendix
\section{Data on degenerations}\label{sec:appendix}
In this appendix, we give explicit transformations which are used in the calculations of degenerations.

\noindent
{\bf 2+1+1 $\to$ 3/2+1+1}

\noindent
$(2)(2), 31, 1111 \to (2)_2, 31, 1111$
\begin{align*}
&\theta^1=2\varepsilon^{-1}, \ \theta^\infty_i=\tilde{\theta}^\infty_i-\varepsilon^{-1} \ (i=1,\ldots,4), \\
&q_1=\tilde{p}_2, \ p_1=\tilde{q}_1-\tilde{q}_2, \ q_2=\tilde{p}_1+\tilde{p}_2, \ p_2=-\tilde{q}_1, \ t=\varepsilon \tilde{t},
\ H=\varepsilon^{-1}\tilde{H}, \\
&Y=(x-1)^{\varepsilon^{-1}}\tilde{Y}.
\end{align*}

\noindent
$(111)(1), 22, 22 \to (1)_2 11, 22, 22$
\begin{align*}
&\theta^\infty_1=\varepsilon^{-1}, \ \theta^\infty_4=\tilde{\theta}^\infty_1-\varepsilon^{-1}, \\
&q_1=\tilde{p}_2, \ p_1=\tilde{q}_1-\tilde{q}_2, \ q_2=\tilde{p}_1+\tilde{p}_2, \ p_2=-\tilde{q}_1,
\ t=\varepsilon \tilde{t}, \ H=\varepsilon^{-1}\tilde{H}, \\
&Y=\tilde{t}^{\varepsilon^{-1}}U^{-1}
\begin{pmatrix} I_2 & O \\
\hat{B}_0 & I_2
\end{pmatrix}
\begin{pmatrix}
0 & 1 & 0 & 0 \\
0 & 0 & q_2 & 0 \\
0 & 0 & 0 & q_1 \\
q_2 & 0 & 0 & 0
\end{pmatrix}\tilde{Y}.
\end{align*}


\noindent
{\bf 3/2+1+1 $\to$ 2+3/2}
\nopagebreak

\noindent
$(2)_2, 31, 1111 \to (111)(1), (2)_2$
\begin{align*}
&\theta^0=-\varepsilon^{-1}, \ \theta^\infty_1=\tilde{\theta}^\infty_1+\varepsilon^{-1}, \\
&q_i=\varepsilon\tilde{t}\tilde{p}_i, \ p_i=-\frac{\tilde{q}_i}{\varepsilon\tilde{t}}\ \  (i=1,2), \ t=-\varepsilon\tilde{t}, \ 
H=-\varepsilon^{-1}\tilde{H}+\frac{\tilde{p}_1\tilde{q}_1+\tilde{p}_2\tilde{q}_2}{\varepsilon\tilde{t}}, \\
&x=1-\varepsilon\tilde{t}\tilde{x}, \ 
Y=\tilde{t}^{-\theta^\infty_1}U^{-1}P^{-1}
{\small
\begin{pmatrix}
1 & -f_4/\theta^0 & \frac{t}{\theta^0}(f_1(1-p_1)+p_1-p_2) & \frac{t}{\theta^0}(1-p_1) \\
0 & 1 & 0 & 0 \\
0 & 0 & t & 0 \\
0 & 0 & 0 & t
\end{pmatrix}
\begin{pmatrix}
1 & 0 & 0 & 0 \\
0 & \tilde{u} & 0 & 0 \\
0 & 0 & \tilde{v} & 0 \\
0 & 0 & 0 & \tilde{w}
\end{pmatrix}}
\tilde{Y}
\end{align*}
where $\tilde{u}, \tilde{v}, \tilde{w}$ satisfiy (\ref{eq:gauge_(1)_2,(2)(2)}).

\noindent
$(1)_2 11, 22, 22 \to (1)_2 11, (2)(2)$
\begin{align*}
&\theta^0=-\varepsilon^{-1}, \ \theta^1=\tilde{\theta}^0+\varepsilon^{-1}, \\
&q_1=\varepsilon\tilde{q}_1,\ p_1=\varepsilon^{-1}\tilde{p}_1, \ q_2=\varepsilon\tilde{q}_2-\frac{\varepsilon\tilde{\theta}^0+1}{\tilde{p}_2}, \ 
p_2=\varepsilon^{-1}\tilde{p}_2, \ t=\varepsilon\tilde{t}, \ H=\varepsilon^{-1}\tilde{H}, \\
&x=\varepsilon^{-1}\tilde{x}, \ 
Y=
\begin{pmatrix}
1 & 0 & 0 & 0 \\
0 & 1 & 0 & 0 \\
0 & 0 & \frac{1}{p_1} & 0 \\
0 & 0 & 0 & -1
\end{pmatrix}\tilde{Y}.
\end{align*}

\noindent
{\bf 2+2 $\to$ 2+3/2}

\noindent
$(111)(1), (2)(2) \to (111)(1), (2)_2$
\begin{align*}
&\theta^0=-2\varepsilon^{-1}, \ \theta^\infty_i=\tilde{\theta}^\infty_i+\varepsilon^{-1} \ (i=1,\ldots,4), \\
&q_i=\varepsilon \tilde{q}_i, \ p_i=\varepsilon^{-1} \tilde{p}_i \ (i=1,2), \ t=\varepsilon \tilde{t}, \ H=\varepsilon^{-1} \tilde{H}, \\
&x=-\varepsilon^{-1}\tilde{x}, \ 
Y=\tilde{x}^{-\varepsilon^{-1}}
\begin{pmatrix}
I_2 & O_2  \\
O_2 & \varepsilon^{-1}I_2
\end{pmatrix}\tilde{Y}.
\end{align*}

\noindent
$(111)(1), (2)(2) \to (1)_2 11, (2)(2)$
\begin{align*}
&\theta^\infty_1=\tilde{\theta}^\infty_1+\varepsilon^{-1}, \ \theta^\infty_3=-\varepsilon^{-1}, \ 
\theta^\infty_4=\tilde{\theta}^\infty_3, \\
&q_1=\varepsilon\left( \tilde{q}_1-\frac{\theta^0+\tilde{\theta}^\infty_1}{\tilde{p}_1} \right), \ p_1=\varepsilon^{-1}\tilde{p}_1, \ 
q_2=\varepsilon\tilde{q}_2-\frac{1}{\tilde{p}_2}, \ p_2=\varepsilon^{-1}\tilde{p}_2, \ t=\varepsilon\tilde{t}, \ 
H=\varepsilon^{-1}\tilde{H}, \\
&Y=U^{-1}
\begin{pmatrix}
I_2 & O_2 \\
B^{(-1)}_0 & I_2
\end{pmatrix}
\begin{pmatrix}
0 & 1 & 0 & 0 \\
0 & 0 & 1 & 0 \\
\varepsilon p_2 & 0 & 0 & 0 \\
0 & 0 & 0 & \varepsilon p_1
\end{pmatrix}\tilde{Y}.
\end{align*}

\noindent
{\bf 3/2+1+1 $\to$ 4/3+1+1}

\nobreak
\noindent
$(1)_2 11, 22, 22 \to (1)_3 1, 22, 22$
\begin{align*}
&\theta^\infty_1=\varepsilon^{-1}, \ \theta^\infty_3=\tilde{\theta}^\infty_1-\varepsilon^{-1}, \\
&q_1=\varepsilon \tilde{t} \tilde{p}_2, \ p_1=-\frac{\tilde{q}_2}{\varepsilon \tilde{t}}, \ 
q_2=\frac{\tilde{t} \tilde{p}_1}{\tilde{p}_1 \tilde{q}_1-\theta^1}, \ 
p_2=-\frac{\tilde{q}_1}{\tilde{t}}(\tilde{p}_1 \tilde{q}_1-\theta^1), \ 
t=-\varepsilon \tilde{t}, \ H=-\varepsilon^{-1}\tilde{H}+\frac{\tilde{p}_1\tilde{q}_1+\tilde{p}_2\tilde{q}_2}{\varepsilon \tilde{t}}, \\
&Y=\frac{1}{p_1+p_2}
\begin{pmatrix}
\varepsilon p_1 & 0 & 0 & -p_1-p_2+\varepsilon \theta^0 p_1 \\
0 & \varepsilon p_1/t & 0 & 0 \\
0 & 0 & -\varepsilon p_1 & \varepsilon p_1 \\
0 & 0 & 0 & -p_1-p_2
\end{pmatrix}\tilde{Y}.
\end{align*}

\noindent
{\bf 3/2+1+1 $\to$ 3/2+1+1}

\noindent
$(1)_2(1)_2, 31, 22 \to ((1))_4, 31, 22$
\begin{align*}
&\theta^\infty_1=\tilde{\theta}^\infty_1+\varepsilon^{-1},\ \theta^\infty_2=-\varepsilon^{-1},\ 
t_1=-\tilde{t}_1,\ t_2=-\tilde{t}_1+\varepsilon \tilde{t}_2,\\
&H_{t_1}=-\tilde{H}_1-\varepsilon^{-1}\tilde{H}_2+\frac{\tilde{p}_2\tilde{q}_2}{\varepsilon \tilde{t}_2},\ 
H_{t_2}=\varepsilon^{-1}\left( \tilde{H}_2-\frac{\tilde{p}_2\tilde{q}_2}{\tilde{t}_2} \right),\\
&q_1=-\tilde{q}_1+\varepsilon \tilde{t}_2 \tilde{p}_2,\ 
p_1=-\frac{\tilde{q}_2}{\varepsilon \tilde{t}_2},\ 
q_2=-\tilde{q}_1,\ p_2=1-\tilde{p}_1+\frac{\tilde{q}_2}{\varepsilon \tilde{t}_2},\\
&Y=f
\begin{pmatrix}
\varepsilon & 1 & 1+\varepsilon(q_1-p_1q_1-p_2q_2-\theta^0-\theta^1) & q_1-q_2(p_1+p_2)-\theta^1 \\
0 & -1 & -1 & (p_1+p_2-1)q_2+\theta^1 \\
0 & 0 & -\varepsilon & -1 \\
0 & 0 & 0 & 1
\end{pmatrix}
\tilde{Y}.
\end{align*}
Here $f$ satisfies
\[
\frac{1}{f}\frac{\partial f}{\partial \tilde{t}_1}=\frac{1}{\tilde{t}_1}\left( 2\tilde{p}_1\tilde{q}_1-\tilde{\theta}^\infty_1+\frac{\tilde{t}_2}{\tilde{t}_1} \right),\quad
\frac{1}{f}\frac{\partial f}{\partial \tilde{t}_2}=-\frac{1}{\tilde{t}_2}\left( 2\tilde{p}_2\tilde{q}_2-\theta^0-\theta^1+\frac{\tilde{t}_2}{\tilde{t}_1} \right).
\]

\noindent
{\bf 3/2+1+1 $\to$ 5/2+1}

\noindent
$(1)_2(1)_2, 31, 22 \to (((1)(1)))_2, 31$
\begin{align*}
&\theta^1=-2\varepsilon^{-3}, \ \theta^\infty_i=\tilde{\theta}^\infty_i+2\varepsilon^{-3},\ 
t_i=-\varepsilon^{-4}\tilde{t}_i-\varepsilon^{-6},\ 
H_{t_i}=-\varepsilon^4 \tilde{H}_i, \\
&q_i=\varepsilon^{-2}\tilde{q}_i-\varepsilon^{-3},\ 
p_i=\varepsilon^2 \tilde{p}_i \ (i=1,2),\\
&x=\varepsilon^2\tilde{x},\ 
Y=\exp[ \varepsilon^{-1}(\tilde{x}-\tilde{t}_1-\tilde{t}_2) ] (\tilde{t}_2-\tilde{t}_1)^{-2\varepsilon^{-3}}
\begin{pmatrix}
1 & 0 & -\varepsilon^{-1} & 0 \\
0 & 1 & 0 & -\varepsilon^{-1} \\
0 & 0 & \varepsilon^2 & 0 \\
0 & 0 & 0 & \varepsilon^2
\end{pmatrix}
\tilde{Y}.
\end{align*}

\noindent
$((1))_4, 31, 22 \to ((((((1))))))_4, 31$
\begin{align*}
&\theta^1=-2\varepsilon^{-3},\ \theta^\infty_1=\tilde{\theta}^\infty_1+4\varepsilon^{-3},\ 
t_1=\varepsilon^{-4}\tilde{t}_1+\varepsilon^{-6},\ t_2=-\varepsilon^{-4}\tilde{t}_2,\\
&H_{t_1}=\varepsilon^4 \tilde{H}_1,\ H_{t_2}=\varepsilon^4 \left( -\tilde{H}_2+\frac{\tilde{p}_2\tilde{q}_2}{\tilde{t}_2} \right),\\
&q_1=\varepsilon^{-3}-\varepsilon^{-2}\tilde{q}_1,\ p_1=1-\varepsilon^2 \tilde{p}_1,\ 
q_2=-\frac{\tilde{t}_2\tilde{p}_2}{\varepsilon^2},\ p_2=\frac{\varepsilon^2\tilde{q}_2}{\tilde{t}_2},\\
&x=\varepsilon^2 \tilde{x},\ 
Y=
\begin{pmatrix}
1 & 0 & 0 & 0 \\
0 & -1 & 0 & 0 \\
0 & 0 & -\varepsilon^2 & 0 \\
0 & 0 & 0 & \varepsilon^2
\end{pmatrix}\tilde{Y}.
\end{align*}

\noindent
\textbf{5/2+1 $\to$ 5/2+1}

\noindent
$(((1)(1)))_2, 31 \to ((((((1))))))_4, 31$
\begin{align*}
&\theta^\infty_1=\tilde{\theta}^\infty_1+\varepsilon^{-1},\ \theta^\infty_2=-\varepsilon^{-1},\ 
t_1=\tilde{t}_1,\ t_2=\tilde{t}_1+\varepsilon \tilde{t}_2,\ 
H_{t_1}=\tilde{H}_1-\varepsilon^{-1}\tilde{H}_2,\ H_{t_2}=\varepsilon^{-1}\tilde{H}_2,\\
&q_1=\tilde{q}_1,\ p_1=\tilde{p}_1-\varepsilon^{-1}\tilde{p}_2,\ q_2=\tilde{q}_1+\varepsilon\tilde{q}_2,\ p_2=\varepsilon^{-1}\tilde{p}_2,\\
&Y=f
\begin{pmatrix}
0 & 1 & 0 & -q_1 \\
-\varepsilon & -1 & \varepsilon q_1 & q_2 \\
0 & 0 & 0 & 1 \\
0 & 0 & -\varepsilon & -1
\end{pmatrix}
\tilde{Y}.
\end{align*}
Here $f$ satisfies
\[
\frac{1}{f}\frac{\partial f}{\partial \tilde{t}_1}=-2\tilde{q}_1,\quad
\frac{1}{f}\frac{\partial f}{\partial \tilde{t}_2}=\frac{1}{\tilde{t}_2}\left( 2\tilde{p}_2\tilde{q}_2+\theta^0 \right).
\]

\noindent
{\bf 2+3/2 $\to$ 3/2+3/2}

\noindent
$(111)(1), (2)_2 \to (1)_2 11, (2)_2$
\begin{align*}
&\theta^\infty_1=\tilde{\theta}^\infty_1-\varepsilon^{-1}, \ \theta^\infty_2=\varepsilon^{-1}, \ 
\theta^\infty_3=\tilde{\theta}^\infty_2, \ \theta^\infty_4=\tilde{\theta}^\infty_3, \ t=\varepsilon\tilde{t}, \ 
H=\varepsilon^{-1}\tilde{H}-\frac{\tilde{p}_2\tilde{q}_2}{\varepsilon\tilde{t}}, \\
&q_1=-\varepsilon \tilde{t} \tilde{p}_2
+\frac{(\tilde{\theta}^\infty_2+\tilde{\theta}^\infty_3)\varepsilon\tilde{t}\tilde{q}_1}{\tilde{q}_1\tilde{q}_2-\tilde{t}}, \ 
p_1=\frac{\tilde{q}_2}{\varepsilon\tilde{t}}, \\
&q_2=-\tilde{q}_1-\frac{2\tilde{\theta}^\infty_2(\tilde{q}_1\tilde{q}_2-\tilde{t})}
{\tilde{p}_1(\tilde{q}_1\tilde{q}_2-\tilde{t})-(\tilde{\theta}^\infty_2+\tilde{\theta}^\infty_3)\tilde{q}_2}
-\frac{(\tilde{q}_1\tilde{q}_2-\tilde{t})^2}{(\tilde{p}_1(\tilde{q}_1\tilde{q}_2-\tilde{t})-(\tilde{\theta}^\infty_2+\tilde{\theta}^\infty_3)\tilde{q}_2)^2}, \ 
p_2=-\tilde{p}_1+\frac{(\tilde{\theta}^\infty_2+\tilde{\theta}^\infty_3)\tilde{q}_2}{\tilde{q}_1\tilde{q}_2-\tilde{t}}, \\
&Y=\tilde{t}^{-\varepsilon^{-1}}U^{-1}
\begin{pmatrix}
I_2 & O_2 \\
B_2 & I_2
\end{pmatrix}\times \\
&{\scriptsize
\begin{pmatrix}
1 & -\frac{\varepsilon}{t}(p_2q_2-2\theta^\infty_3+\frac{1}{p_2})+\frac{1}{t} &
\varepsilon p_1\left( p_2q_2-2\theta^\infty_3+\frac{1}{p_2} \right)+p_2 & 0 \\
0 & 1/t & 0 & 0 \\
-p_2 & 0 & 0 & 0 \\
\frac{p_2(p_1q_2+\varepsilon^{-1})}{\theta^\infty_3-\theta^\infty_4}-p_1 & 0 & 0 &
\frac{p_1}{p_2}\left( p_2q_2-2\theta^\infty_3+\frac{1}{p_2} \right)+\varepsilon^{-1}
\end{pmatrix}}
\tilde{Y}.
\end{align*}

\noindent
$(1)_2 11, (2)(2) \to (1)_2 11, (2)_2$
\begin{align*}
&\theta^0=-2\varepsilon^{-1}, \ \theta^\infty_1=\tilde{\theta}^\infty_1+2\varepsilon^{-1}, \ 
\theta^\infty_2=\tilde{\theta}^\infty_2+\varepsilon^{-1}, \ \theta^\infty_3=\tilde{\theta}^\infty_3+\varepsilon^{-1}, \ 
t=-\varepsilon \tilde{t}, \ H=-\varepsilon^{-1}\tilde{H}+\frac{\tilde{p}_2\tilde{q}_2}{\varepsilon\tilde{t}}, \\
&q_1=\varepsilon \tilde{t} \tilde{p}_2, \ p_1=-\frac{\tilde{q}_2}{\varepsilon \tilde{t}}+\frac{1}{\varepsilon\tilde{q}_1}, \ 
q_2=\tilde{q}_1, \ 
p_2=\tilde{p}_1+\frac{\tilde{\theta}^\infty_1+\tilde{\theta}^\infty_2-\varepsilon^{-1}}{\tilde{q}_1}-\frac{\tilde{t}\tilde{p}_2}{\tilde{q}_1^2},\\
&x=-\varepsilon^{-1}\tilde{x}, \ 
Y=\tilde{t}^{\varepsilon^{-1}}\tilde{x}^{-\varepsilon^{-1}}
\begin{pmatrix}
1 & 0 & 0 & 0 \\
0 & -\varepsilon/t & 0 & 0 \\
0 & 0 & -\varepsilon p_1 & 0 \\
0 & 0 & 0 & q_2
\end{pmatrix}\tilde{Y}.
\end{align*}

\noindent
{\bf 2+3/2 $\to$ 2+4/3}
\nopagebreak

\noindent
$(1)_2 11, (2)(2) \to (1)_3 1, (2)(2)$
\begin{align*}
&\theta^\infty_1=\varepsilon^{-1}, \ \theta^\infty_3=\tilde{\theta}^\infty_1-\varepsilon^{-1}, \\
&q_1=\varepsilon\tilde{t}\tilde{p}_1, \ p_1=-\frac{\tilde{q}_1}{\varepsilon\tilde{t}}, \ t=-\varepsilon\tilde{t}, \ 
H=-\varepsilon^{-1}\tilde{H}+\frac{\tilde{p}_1\tilde{q}_1}{\varepsilon\tilde{t}}, \\
&Y=\tilde{t}^{\varepsilon^{-1}}
\begin{pmatrix}
1 & 0 & 0 & 1/t \\
0 & 1/t & 0 & 0 \\
0 & 0 & p_1 & 0 \\
0 & 0 & 0 & -1/t
\end{pmatrix}
\tilde{Y}.
\end{align*}

\noindent
{\bf 4/3+1+1 $\to$ 2+4/3}

\noindent
$(1)_3 1, 22, 22 \to (1)_3 1, (2)(2)$
\begin{align*}
&\theta^0=-\varepsilon^{-1}, \ \theta^1=\tilde{\theta}^0+\varepsilon^{-1}, \\
&q_1=\frac{\tilde{t}}{\tilde{q}_2}, \ p_1=-\frac{\tilde{q}_2}{\tilde{t}}(\tilde{p}_2\tilde{q}_2-\tilde{\theta}^0-\varepsilon^{-1}), \ 
q_2=\tilde{q}_1, \ p_2=\tilde{p}_1, \ t=\varepsilon\tilde{t}, \ H=\varepsilon^{-1}\tilde{H}-\frac{\tilde{p}_2\tilde{q}_2}{\varepsilon\tilde{t}}, \\
&x=\varepsilon^{-1}\tilde{x}, \ 
Y=
\begin{pmatrix}
\varepsilon^{-1} & -p_2-\frac{\theta^0}{q_1} & \theta^0\frac{q_2}{q_1}+\theta^1+\theta^\infty_2 & -\theta^0/t \\
0 & 1 & 0 & 0 \\
0 & 1/q_1 & -q_2/q_1 & 1/t \\
0 & 1/q_1 & 1-\frac{q_2}{q_1} & 1/t
\end{pmatrix}
\tilde{Y}.
\end{align*}

\noindent
{\bf 4/3+1+1 $\to$ 5/4+1+1}
\nopagebreak

\noindent
$(1)_3 1, 22, 22 \to (1)_4, 22, 22$
\begin{align*}
&\theta^\infty_1=\varepsilon^{-1}, \ \theta^\infty_2=\tilde{\theta}^\infty_1-\varepsilon^{-1}, \ 
t=-\varepsilon \tilde{t}, \ H=-\varepsilon^{-1}\tilde{H}+\frac{\tilde{p}_2\tilde{q}_2}{\varepsilon\tilde{t}}, \\
&q_1=\varepsilon \tilde{t} \tilde{p}_2, \ p_1=-\frac{\tilde{q}_2}{\varepsilon \tilde{t}}+\frac{1}{\varepsilon\tilde{q}_1}, \ 
q_2=\tilde{q}_1, \ 
p_2=\tilde{p}_1+\frac{\theta^0+\theta^1+\tilde{\theta}^\infty_1-\varepsilon^{-1}}{\tilde{q}_1}-\frac{\tilde{t}\tilde{p}_2}{\tilde{q}_1^2},\\
&Y=
\begin{pmatrix}
1 & -p_2+\frac{p_2q_1+\theta^0+\theta^1+\theta^\infty_2}{q_2}-\frac{\theta^1 p_2}{p_1q_2}-\frac{\theta^1}{\varepsilon p_1 q_2^2} &
\frac{1}{t}\left( \varepsilon(p_2q_1+\theta^0)-\frac{\theta^1}{p_1q_2}-\frac{\varepsilon \theta^1 p_2}{p_1} \right) & 0 \\
0 & 1-\frac{q_1}{q_2} & -\frac{\varepsilon q_1}{t} & -\varepsilon \\
0 & -\frac{1}{q_2} & -\frac{\varepsilon}{t} & 0 \\
0 & 0 & -\frac{\varepsilon}{t} & 0
\end{pmatrix}
\tilde{Y}.
\end{align*}

\noindent
{\bf 3/2+3/2 $\to$ 3/2+4/3}
\nopagebreak

\noindent
$(1)_2 11, (2)_2 \to (1)_3 1, (2)_2$
\begin{align*}
&\theta^\infty_1=-\varepsilon^{-1}, \ \theta^\infty_2=\tilde{\theta}^\infty_2+\frac12, \ 
\theta^\infty_3=\tilde{\theta}^\infty_1-\frac12+\varepsilon^{-1}, \\
&q_1=\tilde{q}_1+\varepsilon\tilde{q}_1\left( \tilde{p}_1\tilde{q}_1-2\tilde{\theta}^\infty_2-\frac12 \right), \ 
p_1=\frac{1}{\varepsilon\tilde{q}_1}, \ 
t=\varepsilon \tilde{t}, \ H=\varepsilon^{-1}\tilde{H}, \\
&Y=\tilde{t}^{-\varepsilon^{-1}+\tilde{\theta}^\infty_2}
\begin{pmatrix}
1 & 0 & 0 & q_1/t \\
0 & 1 & 0 & 0 \\
0 & 0 & 1 & 0 \\
0 & 0 & 0 & -1/t
\end{pmatrix}\tilde{Y}.
\end{align*}

\noindent
{\bf 2+4/3 $\to$ 3/2+4/3}

\noindent
$(1)_3 1, (2)(2) \to (1)_3 1, (2)_2$
\begin{align*}
&\theta^0=2\tilde{\theta}^\infty_2+2\varepsilon^{-1}, \ 
\theta^\infty_1=\tilde{\theta}^\infty_1-3\tilde{\theta}^\infty_2-\frac12-3\varepsilon^{-1}, \ 
\theta^\infty_2=\frac12-\varepsilon^{-1}, \\
&q_1=\tilde{q}_2, \ p_1=\tilde{p}_2, \ q_2=\tilde{q}_1, \ p_2=\tilde{p}_1+\frac{1}{\varepsilon\tilde{q}_1}, \ 
t=\varepsilon\tilde{t}, \ H=\varepsilon^{-1}\tilde{H}, \\
&x=\varepsilon^{-1}\tilde{x}, \ 
Y=\tilde{x}^{\tilde{\theta}^\infty_2+\varepsilon^{-1}}\tilde{t}^{\,2\tilde{\theta}^\infty_2+\varepsilon^{-1}}
\begin{pmatrix}
\varepsilon^{-1} & 0 & 0 & 0 \\
0 & 1 & 0 & 0 \\
0 & 0 & 1 & 0 \\
0 & 0 & 0 & q_2
\end{pmatrix}\tilde{Y}.
\end{align*}

\noindent
{\bf 2+4/3 $\to$ 2+5/4}

\noindent
$(1)_3 1, (2)(2) \to (1)_4, (2)(2)$
\begin{align*}
&\theta^\infty_1=\tilde{\theta}^\infty_1-1+\varepsilon^{-1}, \ \theta^\infty_2=1-\varepsilon^{-1}, \\
&q_1=-\tilde{q}_1\left( \varepsilon \tilde{p}_1\tilde{q}_1+\frac{\varepsilon}{2}-1 \right), \ 
p_1=-\frac{1}{\varepsilon\tilde{q}_1}, \ t=-\varepsilon\tilde{t}, \ H=-\varepsilon^{-1}\tilde{H}, \\
&Y=
\begin{pmatrix}
1 & 0 & 0 & 0 \\
0 & tp_1q_1 & t & -tp_2 \\
0 & tp_1 & 0 & 0 \\
0 & 0 & 0 & t
\end{pmatrix}\tilde{Y}.
\end{align*}

\noindent
{\bf 5/4+1+1 $\to$ 2+5/4}

\noindent
$(1)_4, 22, 22 \to (1)_4, (2)(2)$
\begin{align*}
&\theta^0=\tilde{\theta}^0-\varepsilon^{-1}, \ \theta^1=\varepsilon^{-1}, \\
&q_1=\tilde{q}_1, \ p_1=\tilde{p}_1-\frac{\tilde{\theta}^0+\frac12}{\tilde{q}_1}-\frac{\tilde{t}}{\tilde{q}_1^2\tilde{q}_2}, \ 
q_2=\tilde{q}_2(1-\varepsilon\tilde{p}_2\tilde{q}_2)+\frac{\varepsilon\tilde{t}}{\tilde{q}_1}, \ p_2=-\frac{1}{\varepsilon\tilde{q}_2}, \ 
t=\varepsilon\tilde{t}, \ H=\varepsilon^{-1}\tilde{H}-\frac{1}{\varepsilon\tilde{q}_1\tilde{q}_2}, \\
&x=\varepsilon^{-1}\tilde{x}, \ 
Y=
\begin{pmatrix}
1 & 0 & 0 & -t p_2/q_1 \\
0 & t & 0 & 0 \\
0 & t^2/q_1 & -t/p_2 & -\frac{\varepsilon t p_2}{q_1}(q_1q_2-t) \\
0 & 0 & 0 & t p_2
\end{pmatrix}\tilde{Y}.
\end{align*}

\noindent
{\bf 3/2+4/3 $\to$ 3/2+5/4}

\noindent
$(1)_3 1, (2)_2 \to (1)_4, (2)_2$
\begin{align*}
&\theta^\infty_1=\frac12-\varepsilon^{-1}, \ \theta^\infty_2=-\frac12+\varepsilon^{-1}, \\
&q_i=-\frac{\tilde{t}}{\tilde{q}_i}, \ p_i=\frac{\tilde{q}_i}{\tilde{t}}\left( \tilde{p}_i\tilde{q}_i-\frac{1}{\varepsilon}+\frac12 \right) \ \ (i=1,2), \ 
t=\varepsilon\tilde{t}, \ H=\varepsilon^{-1}\tilde{H}-\frac{\tilde{p}_1\tilde{q}_1+\tilde{p}_2\tilde{q}_2}{\varepsilon\tilde{t}}, \\
&Y=\tilde{t}^{-\varepsilon^{-1}+\frac12}
\begin{pmatrix}
1 & 0 & 0 & 0 \\
0 & t/\varepsilon & 0 & 0 \\
0 & \frac{t}{\varepsilon q_2} & \frac{t}{q_1q_2} & 0 \\
0 & 0 & 0 & -t/q_1
\end{pmatrix}\tilde{Y}.
\end{align*}

\noindent
{\bf 2+5/4 $\to$ 3/2+5/4}

\noindent
$(1)_4, (2)(2) \to (1)_4, (2)_2$
\begin{align*}
&\theta^0=2\varepsilon^{-1}, \ \theta^\infty_1=-4\varepsilon^{-1}, \\
&q_1=-\frac{\tilde{t}}{\tilde{q}_2}, \ p_1=\frac{\tilde{q}_2}{\tilde{t}}\left( \tilde{p}_2\tilde{q}_2-\varepsilon^{-1} \right), \ 
q_2=-\frac{\tilde{t}}{\tilde{q}_1}, \ p_2=\frac{\tilde{q}_1}{\tilde{t}}\left( \tilde{p}_1\tilde{q}_1-\varepsilon^{-1} \right), \ 
t=\varepsilon\tilde{t}, \ H=\varepsilon^{-1}\tilde{H}-\frac{\tilde{p}_1\tilde{q}_1+\tilde{p}_2\tilde{q}_2}{\varepsilon\tilde{t}}, \\
&x=\varepsilon^{-1}\tilde{x}, \ 
Y=
\tilde{x}^{\varepsilon^{-1}}\tilde{t}^{\varepsilon^{-1}}
\begin{pmatrix}
1 & 0 & 0 & 0 \\
0 & 1 & 0 & 0 \\
0 & 0 & -1/q_2 & -p_2 \\
0 & 0 & 0 & -1
\end{pmatrix}\tilde{Y}.
\end{align*}

\end{document}